\date{}

\documentclass[11pt]{article}
\usepackage{amsfonts}
\usepackage{graphicx}
\usepackage{latexsym}
\usepackage{amssymb}
\usepackage{amsmath}
\usepackage{enumerate}
\usepackage{layout}
\usepackage{dsfont}
\usepackage{eufrak}
\usepackage{color}
\usepackage{hyperref}
\usepackage{cite}
\usepackage{color}
\usepackage[normalem]{ ulem }
\usepackage{soul}

\newtheorem{prop}{Proposition}
\newtheorem{lemma}{Lemma}

\newtheorem{theorem}{Theorem}
\newtheorem{remark}{Remark}

\def\real{{\mathord{{\rm I\kern-2.8pt R}}}}
\def\inte{{\mathord{{\rm I\kern-2.8pt N}}}}

\def\sZZ{{\rm Z\kern-2.8ptem{}Z}}

\def\z{{\mathchoice
  {\sZZ}
  {\sZZ}
  {\rm Z\kern-0.30em{}Z}
  {\rm Z\kern-0.25em{}Z} }}
\def\sQQ{{\kern 0.27em \vrule height1.45ex width0.03em depth0em
          \kern-0.30em \rm Q}}
\def\qu{{\mathchoice
    {\sQQ}
    {\sQQ}
  {\kern 0.225em \vrule height1.05ex width0.025em depth0em \kern-0.25em \rm Q}
  {\kern 0.180em \vrule height0.78ex width0.020em depth0em \kern-0.20em \rm Q}
        }}
\def\sCC{{\kern 0.27em \vrule height1.45ex width0.03em depth0em
          \kern-0.30em \rm C}}
\def\complex{{\mathchoice
    {\sCC}
    {\sCC}
  {\kern 0.225em \vrule height1.05ex width0.025em depth0em \kern-0.25em \rm C}
  {\kern 0.180em \vrule height0.78ex width0.020em depth0em \kern-0.20em \rm C}
        }}

\newcommand{\ba}{\begin{array}}
\newcommand{\ea}{\end{array}}
\newcommand{\be}{\begin{equation}}
\newcommand{\ee}{\end{equation}}
\newcommand{\bea}{\begin{eqnarray}}
\newcommand{\eea}{\end{eqnarray}}
\newcommand{\beaa}{\begin{eqnarray*}}
\newcommand{\eeaa}{\end{eqnarray*}}

\newcommand{\eps}{\varepsilon}

\def\b{\beta}

\def\z{\zeta}

\font\tenmath=msbm10 \font\sevenmath=msbm7 \font\fivemath=msbm5
\newfam\mathfam \textfont\mathfam=\tenmath
\scriptfont\mathfam=\sevenmath \scriptscriptfont\mathfam=\fivemath

\def \b{\noindent}

\def \={{\buildrel {\rm (law)} \over =}}

\def\qed{ \hfill \vrule width.25cm height.25cm depth0cm\smallskip}

\newcommand{\basa}{\begin{assumption}}
\newcommand{\easa}{\end{assumption}}

\newcommand{\bas}{\begin{assum}}
\newcommand{\eas}{\end{assum}}

\def\span{\hbox{\rm span$\,$}}

\newcommand{\ignore}[1]{}
\allowdisplaybreaks
\begin{document}

\title{ Behavior  with respect to the Hurst index of the Wiener Hermite  integrals and application to SPDEs  }
\author{Meryem Slaoui and C. A. Tudor \vspace*{0.2in} \\
 $^{1}$ Laboratoire Paul Painlev\'e, Universit\'e de Lille 1\\
 F-59655 Villeneuve d'Ascq, France.\\
\quad meryem.slaoui@univ-lille.fr\\
 \quad ciprian.tudor@univ-lille.fr\\
\vspace*{0.1in} }

\maketitle

\begin{abstract}
We consider the Wiener integral with respect to a $d$-parameter Hermite process with Hurst multi-index ${\bf H}= (H_{1},.., H_{d}) \in \left( \frac{1}{2}, 1\right) ^{d}$ and we analyze the limit behavior in distribution of this object when the components of ${\bf H}$ tend to $1$ and/or $\frac{1}{2}$. As examples, we focus on the solution to the stochastic heat equation with additive Hermite noise and to the Hermite Ornstein-Uhlenbeck process.
\end{abstract}

\vskip0.3cm
{\bf 2010 AMS Classification Numbers:}   60H05, 60H15, 60G22.

\vskip0.3cm


\section{Introduction}

The Hermite processes are self-similar processes with long-memory and stationary increments. These properties made them good models for many applications. The Hermite processes constitute a non-Gaussian extension of the fractional Brownian motion. Their Hurst parameter, which is contained in the interval $\left(\frac{1}{2}, 1\right)$, characterizes the main properties of this process. The reader may consult the monographs \cite{PiTa-book} or \cite{T} for a complete exposition on Hermite processes.

Our work deals with stochastic partial differential equations (SPDEs) driven by the Hermite process. Starting with the seminal work \cite{Walsh1986}, many researchers explored the possibility of solving SPDEs with general  noises more general than the standard space-time white noise. In our work, such a stochastic perturbation is chosen to be  the Hermite noise.    Recently, various types of stochastic integral and stochastic equations driven by Hermite noises have been considered by many authors. We refer, among others, to \cite{Ba}, \cite{Cou}, \cite{CouMas}, \cite{CMO}, \cite{NT}, \cite{Tran}, \cite{BonTud}, \cite{GIP}, \cite{SlaTud1}, \cite{SlaTud2}. Our purpose is to analyze the asymptotic behavior in distribution of the solution to the stochastic heat equation with additive Hermite noise, when the Hurst parameter (which is also the self-similarity index of the Hermite process) converges to the extreme values of its interval of definition, i.e when it tends to one and to one half. Our work continues a recent line of research that concerns the limit behavior in distribution with respect to the Hurst parameter of Hermite and related  fractional-type stochastic processes. In particular, the papers \cite{BellNu} and \cite{BaiTa} deal with the asymptotic behavior of the generalized Rosenblatt process, the work \cite{AT} studies the multiparamter Hermite processes while the paper \cite{SlaTud2} investigates the Ornstein-Uhlenbeck process with Hermite noise of order $q=2$.

The solution to the heat equation with Hermite noise in $\mathbb{R} ^{d}$ is a $(d+1)$- parameter random field depending on a Hurst index ${\bf H} \in \left( \frac{1}{2}, 1\right) ^{d+1}$.  We prove that the solution converges in distribution to a Gaussian limit when at least one of the components of ${\bf H}$ converges to $\frac{1}{2}$ and to a random variable in a Wiener chaos of higher order when at least one of the components of ${\bf H}$ tends to $1$ (and none of them converges to $\frac{1}{2}$). Moreover, the limit always coincides in distribution with the solution to the stochastic heat equation driven by the limit of the Hermite noise. The results show that these models offer a large flexibilitily, covering a large class of probability distributions, from Gaussian laws to distribution of random variables in Wiener chaos of higher order. 

For the proofs  we use various techniques, such as the Malliavin calculus and the Fourth Moment Theorem for the normal convergence, the properties of the Wiener  integrals with respect to the Hermite process and the so-called power counting theorem.  Since the solution to the Hermite-driven heat equation can be expressed as a Wiener integral with respect to a Hermite sheet, we start our analysis by some more general results, i.e by studying the   behavior with respect to the Hurst index of such Wiener integrals. This allows to consider other examples, in particular the Hermite Ornstein-Uhlenbeck process. 

We organized our paper as follows. Section 2 contains some preliminaries. We introduce the multidimensional  Hermite processes  and the Wiener integral with respect to them. We also recall some known results  concerning the asymptotic behavior of the Hermite sheet. In Section 3, we state general results on the asymptotic behavior of the Wiener-Hermite integrals with respect to  the Hurst parameter.  We will give two applications of the main results obtained. In Section 4 we analyse the asymptotic behavior of the mild  solution of the stochastic heat equation with Hermite noise and finally Section 5 contains the case of  the  Hermite Ornstein -Uhlenbeck process. The Appendix (Section 6) contains the basic elements of the stochastic analysis on Wiener spaces needed in the paper.

\section{Preliminaries}
 In this preliminary  section we will introduce  the Hermite sheet and the Wiener integral with respect to this multiparameter process. We also recall the main findings from \cite{AT} concerning the behavior of the Hermite sheet with respect to its Hurst multi-index. We start with some multidimensional notation, that we will use throughout our work.

\subsection{Notation}
For $d\in \mathbb{N} \backslash \left\lbrace0\right\rbrace$ we will work with multi-parametric processes indexed by elements of  $\mathbb{R} ^{d}$. We shall use bold notation for multi-indexed quantities, i.e., $\mathbf{a}=(a_{1},a_{2},\ldots ,a_{d})$, $ \mathbf{b}= (b_{1},b_{2},.., b_{d})$, $\mathbf{\alpha }= (\alpha_{1}, ..,\alpha_{d})$,  $\mathbf{ab}=\prod_{i=1} ^{d} a_{i} b_{i}$, $\vert \mathbf{a}- \mathbf{b} \vert  ^{\mathbf{\alpha }}  =\prod _{i=1} ^{d} \vert a_{1}- b_{1}\vert ^{\alpha _{i}}$, $\mathbf{a/b}=(a_{1}/b_{1},a_{2}/b_{2},\ldots ,a_{d}/b_{d})$, $ [\mathbf{a},\mathbf{b}]= \displaystyle \prod_{i=1}^{d}[a_{i},b_{i}]$, $(\mathbf{a},\mathbf{b})=\displaystyle \prod_{i=1}^{d}(a_{i},b_{i})$, $\displaystyle\sum_{\mathbf{i} =0 }^{\mathbf{N}} a_{\mathbf{i}} = \displaystyle\sum_{i_{1}=0}^{N_{1}} \displaystyle\sum_{i_{2}=0}^{N_{2}} \ldots \displaystyle\sum_{i_{d}=0}^{N_{d}} a_{i_{1},i_{2},\ldots ,i_{d}}$ if ${\bf N}=(N_{1},.., N_{d})$, $\mathbf{a}^{\mathbf{b}}=\displaystyle \prod_{i=1}^{d} a_{i}^{b_{i}}$, and $ \mathbf{a} < \mathbf{b} $ iff $ a_{1} < b_{1}, a_{2} < b_{2},\ldots ,a_{d} < b_{d}$ (analogously for the other inequalities). 

We write $ {\bf a}-{\bf 1}$ to indicate the product $\prod_{i=1} ^{d} (a_{i}-1).$ By $\beta $ we denote the  Beta function $\beta (p,q)= \int_{0}^{1} z^{p-1}(1-z)^{q-1} dz, p,q>0$ and  we use  the notation $$\beta ({\bf a}, {\bf b})= \prod_{i=1}^{d} \beta \left( a^{(i)}, b^{(i)}\right)$$
if ${\bf a}= (a^{(1)},.., a^ {(d) })$ and ${\bf b} = (b^{(1)},.., b^{(d)})$.  

 Let us recall that the  increment of a $d$-parameter process $X$ on a rectangle $[\mathbf{s}, \mathbf{t} ] \subset \mathbb{R} ^{d}$, $\mathbf{s}= (s_{1},\ldots , s_{d} ), \mathbf{t}=(t_{1},\ldots ,t_{d})$, with $\mathbf{s} \leq \mathbf{t}$ (denoted by $\Delta X ([\mathbf{s}, \mathbf{t} ])$) is given by
\begin{equation}\label{mi}
\Delta X ([\mathbf{s}, \mathbf{t} ])= \displaystyle\sum_{r\in \{0,1\} ^{d}}(-1) ^{d- \sum_{i=1}^{d} r_{i}} X_{\mathbf{s} + \mathbf{r} \cdot (\mathbf{t}-\mathbf{s})}.
\end{equation}
When $d=1$ one obtains $\Delta X ([\mathbf{s}, \mathbf{t} ])=X_{t}-X_{s}$ while for $d=2$ one gets $\Delta X ([\mathbf{s}, \mathbf{t} ])=X_{t_{1}, t_{2}} -X_{t_{1}, s_{2}} -X_{s_{1}, t_{2}} + X_{s_{1}, s_{2}}$.
\subsection{Hermite processes and Wiener-Hermite integrals}

We recall the definition and the basic properties of multiparameter Hermite processes. For a more complete presentation, we refer to \cite{DCT}, \cite{PiTa-book} or \cite{T}.
Let $q\geq 1$ integer and the Hurst multi-index $\mathbf{H}=(H_{1},H_{2}, \ldots ,H_{d}) \in (\frac{1}{2},1)^{d}$. The {\em Hermite sheet of order q and with self-similarity index {\bf H} }, denoted $ ( Z ^{q,d} _{{\bf H}} ({\bf t}), {\bf t} \in \mathbb{R} _{+} ^{d})$ in the sequel,  is given by

\begin{eqnarray}
\label{hermite-sheet-1}
\nonumber
 Z^{q, d}_{\mathbf{H}}(\mathbf{t}) &=& c(\mathbf{H},q) \int_{\mathbb{R}^{d\cdot q}}
  \int_{0}^{t^{(1)}}
 \ldots \int_{0}^{t^{(d)}} \left( \prod _{j=1}^{q} (s_{1}-y_{1,j})_{+} ^{-\left( \frac{1}{2} + \frac{1-H_{1}}{q} \right) }
 \ldots (s_{d}-y_{d,j})_{+} ^{-\left( \frac{1}{2} + \frac{1-H_{d}}{q} \right) }   \right)   \\
\nonumber
& &   ds_{d} \ldots ds_{1}  \quad
dW(y_{1,1},\ldots ,y_{d,1} )\ldots dW(y_{1,q},\ldots ,y_{d,q}) \\
&=&
c(\mathbf{H},q) \int_{\mathbb{R}^{d\cdot q}}
\int_{0}^{\mathbf{t}} \quad \prod _{j=1}^{q} (\mathbf{s}-\mathbf{y}_{j})_{+} ^{-\mathbf{\left( \frac{1}{2} + \frac{1-H}{q} \right)} }  d\mathbf{s} \quad dW(\mathbf{y}_{1})\ldots dW(\mathbf{y}_{q})
\end{eqnarray}
for every ${\bf t} =(t^{1},..., t^{d})\in \mathbb{R} _{+} ^{d}$,  where $x_{+}=\max(x,0)$. The above   stochastic integral is  a multiple stochastic integral with respect to the Wiener sheet  ($W(\mathbf{y}), {\bf y} \in \mathbb{R} ^{d}$), see Section \ref{app1}.  The constant $c({\bf H}, q)$ ensures that $ \mathbf{E} \left( Z^{q}_{\mathbf{H}}(\mathbf{t}) \right) ^{2} ={\bf t}^{2{\bf H}}$ for every ${\bf t} \in \mathbb{R} ^{d}_{+}$. As pointed out before, when $q=1$, (\ref{hermite-sheet-1}) is the fractional Brownian sheet with Hurst multi-index $\mathbf{H}=(H_{1},H_{2}, \ldots ,H_{d}) \in (\frac{1}{2},1)^{d}$. For $q \geq 2$ the process $Z^{q,d} _{\mathbf{H}}$ is not Gaussian and for $q=2$ we denominate it as the {\em Rosenblatt sheet}.

The Hermite sheet is a ${\bf H}$-self-similar stochastic process and it has stationary increments. Its paths are H\"older continuous of order $\boldsymbol {\delta} < {\bf H}$, see \cite{PiTa-book} or \cite{T}. Its covariance is the same for every $q\geq1$  and it coincides with the covariance of the $d$-parameter fractional Brownian motion, i.e.

\begin{equation}
\mathbf{E}Z _{{\bf H} }^ {q, d} ({\bf t}) Z _{{\bf H} }^ {q, d} ({\bf s}) =\prod _{j=1} ^ {d} \left( \frac{1}{2} \left( t _{i} ^ {2H_{i} }+s_{i} ^ {2H_{i}}-\vert t_{i}-s_{i} \vert ^ {2H_{i}}\right) \right)=: R_{{\bf H}}({\bf t}, {\bf s}),  \ \ \ \  t_{i},s_{i} \geq 0.\label{covh}
\end{equation}

We will denote by $\vert \mathcal{H} _{{\bf H}}\vert$ the space of measurable functions $f:\mathbb{R} ^{d}\to \mathbb{R} $ such that  
$$\Vert f \Vert ^{2} _{ \vert \mathcal{H} _{{\bf H}}\vert } <\infty$$
where
\begin{eqnarray}
\label{normahb}
\Vert f \Vert ^{2} _{ \vert \mathcal{H} _{{\bf H}}\vert } &: =& {\bf H} (2{\bf H}-{\bf 1}) \int_{\mathbb{R} ^{d} }  \int_{\mathbb{R} ^{d} }  d {\bf u} d{\bf v} \vert f({\bf u})\vert \cdot \vert  f ({\bf v}) \vert \vert {\bf u} -{\bf v}\vert ^{2{\bf H}-2}\\
&=& {\bf H} (2{\bf H}-{\bf 1}) \int_{\mathbb{R} ^{d} }  \int_{\mathbb{R} ^{d} } du^{(1) }...du^{(d)} dv^{(1) }...dv^{(d)}\nonumber\\
&&\times  f(u^{(1)},.., u^{(d)})f(v^{(1)},.., v^{(d)})\prod_{j=1}^{d} \vert u ^{(j)}-v ^{(j)} \vert ^{2H_{j}-2} \nonumber
\end{eqnarray}
where ${\bf u}= (u ^{(1)},.., u^{(d)}), {\bf v}=(v ^{(1)},.., v ^{(d)})\in \mathbb{R} ^{d}.$

Notice that the space $\vert \mathcal{H}_{{\bf H}}\vert $ satisfies the following inclusion (see Remark 3 in \cite{DCT})

\begin{equation}
\label{inclu}
L ^{1} (\mathbb{R} ^{d}) \cap L ^{2} (\mathbb{R} ^{d})  \subset  L ^{\frac{1}{{\bf H}}} (\mathbb{R} ^{d}) \subset \vert \mathcal{H}_{{\bf H}}\vert.
\end{equation}

The Wiener integral with respect to the Hermite sheet $Z ^{q,d}_{{\bf H}}$ has been defined in \cite{DCT} (following the idea of \cite{MaTu} in the one-parmeter case). In  particular,  it is well-defined for measurable integrands  $f\in \vert \mathcal{H}_{{\bf H}}\vert $ via the formula

\begin{equation}\label{hw}
\int_{\mathbb{R} ^{d}} f({\bf s}) dZ ^{q,d}_{{\bf H}} ({\bf s}) = \int_{\mathbb{R} ^{d. q}} (Jf)({\bf y}_{1},..., {\bf y}_{q} )dW({\bf y}_{1})...dW ({\bf y}_{q}) 
\end{equation}
where $ \left(W ({\bf y}), {\bf y}\in \mathbb{R} ^{d} \right)$ is a $d$-parameter Wiener process and 
\begin{equation}
\label{jf}
(Jf)({\bf y}_{1},..., {\bf y}_{q} ) =c({\bf H}, q)  \int_{\mathbb{R} ^ {d}   }d{\bf u} f({\bf u}) ({\bf u}-{\bf y }_{1} )_{+}^ {-\left( \frac{1}{2}+ \frac{1-{\bf H} }{q} \right)} \ldots ({\bf u}-{\bf y }_{q} )_{+}^ {-\left( \frac{1}{2}+ \frac{1-{\bf H} }{q} \right)}
\end{equation}
with $ c({\bf H}, q)  $ from (\ref{hermite-sheet-1}). The stochastic integral $\int_{\mathbb{R} ^{d. q}} (Jf)({\bf y}_{1},..., {\bf y}_{q} )dW({\bf y}_{1})...dW ({\bf y}_{q}) $ is a multiple Wiener-It\^o integral with respect to the Wiener sheet $W$.

We have the isometry formula, for $f,g\in  \vert \mathcal{H} _{{\bf H}}\vert$
\begin{eqnarray}
\mathbf{E} \left(\int_{\mathbb{R} ^{d}} f({\bf s}) dZ ^{q,d}_{{\bf H}} ({\bf s}) \int_{\mathbb{R} ^{d}} g({\bf s}) dZ ^{q,d}_{{\bf H}} ({\bf s}) \right) &=& {\bf H} (2{\bf H}-{\bf 1}) \int_{\mathbb{R} ^{d} }  \int_{\mathbb{R} ^{d} }  d {\bf u} d{\bf v}  f({\bf u})   g ({\bf v})  \vert {\bf u} -{\bf v}\vert ^{2{\bf H}-2}\nonumber \\
& :=& \langle f,g\rangle _{\mathcal {H}_{ {\bf H}} }.\label{isoHW}
\end{eqnarray} 
By $\Vert f\Vert _{\mathcal{H}_{\bf H}}^{2}$ we denote $\langle f,f\rangle _{\mathcal {H}_{ {\bf H}} }$.

\subsection{ Behavior of the Hermite sheet with respect to the Hurst parameter}

In a first step, we analyze the convergence of the integral $ \int_{\mathbb{R} ^{d}} f({\bf s}) dZ ^{q,d}_{{\bf H}} ({\bf s}) $ when the Hurst indices $H_{i}$ goes to 1 and/or $\frac{1}{2}$.

Let us introduce the following notation: if $\{ j_{1},.., j_{k} \} \subset \{ 1,.., d\}$ with $1\leq k\leq d$ we will denote

\begin{equation}
\label{not}
A_{k}= \{ j_{1},.., j_{k}\}, \hskip0.3cm  {\bf H}_{A_{k}}= (H_{j_{1}},..., H_{j_{k}}) \in \left(\frac{1}{2}, 1\right) ^{k} , \hskip0.3cm \langle {\bf t}\rangle _{A_{k}}= t ^{(j_{1})}....t ^{(j_{k})} \mbox{ if } {\bf t} = ( t ^{(1)},.., t^{(d)}).
\end{equation}

We will separate our study into  following two situations:

\begin{enumerate}

\item At least one parameter converges to $1$ and none to $\frac{1}{2}$. Then the limit will be a non-Gaussian random variable related to the Hermite distribution. 

\item At least one parameter $H_{i}$ converges to $\frac{1}{2}$ and the other indices are fixed in $(\frac{1}{2}, 1)$ or converges to 1, i.e. if $A_{k}$ is as above, $B_{p}= \{l_{1},.., l_{p} \}\subset \{1,..,d\} $ with $p+k\leq d$ and $A_{k}\cap B_{p} =\emptyset$, we assume
 ${\bf H} _{A_{k}}\to (\frac{1}{2},..., \frac{1}{2}) \in \mathbb{R} ^{k}$ and ${\bf H}_{B_{p}} \to (1,..,1) \in \mathbb{R} ^{p}. $  In this case we will see that the limit of $ \int_{\mathbb{R} ^{d}} f({\bf s}) dZ ^{q,d}_{{\bf H}} ({\bf s}) $ is a centered Gaussian random variable with an explicit variance.

\end{enumerate}

We start by recalling the main result in \cite{AT} concerning the asymptotic behavior of the Hermite sheet.

\begin{theorem}\label{tat}  Let $\left( Z ^{q, d} _{{\bf H}} ({\bf t}) \right) _{{\bf t}\geq 0}$ be given by  (\ref{hermite-sheet-1}) and let $A_{k}, B_{p} $ be as in (\ref{not}). Fix $T>0$.

\begin{enumerate}
\item 
Assume ${\bf H} _{A_{k}}\to \left( 1,..,1\right) \in \mathbb{R} ^{k}$. Assume that the parameters $H_{j}, j\in \overline{A}_{k}$ are fixed.  Then the process $ Z^{q,d}_ {{\bf H}}$ converges weakly  in $C([0,T] ^{d} )$ to the $d$-parameter stochastic process  $(X_{{\bf t}})_{{\bf t}\geq 0}$ defined by
\begin{equation}\label{30s-2}
X_{{\bf t}}=  \langle {\bf t}\rangle _{A_{k}} Z _{{\bf H} _{\overline{A}_{k}}}^{q, d-k} ( {\bf t}_{\overline{A}_{k}})
\end{equation}
where $\left( Z_{{\bf H}_{\overline{A}_{k}}} ^{q, d-k} ({\bf  t} _{\overline{A}_{k}})\right) _{ {\bf t}_{\overline{A}_{k}}\in \mathbb{R}^{d-k} _{+} }$ is a $(d-k)$-parameter Hermite process of order $q$ with Hurst index $ {\bf H}_{\overline{A}_{k}}\in \left(\frac{1}{2}, 1\right) ^{d-k}$.

\item Assume $(H_{1},.., H_{d}) \to (1,..,1) \in \mathbb{R} ^{d}$. Then the process $ Z^{q, d }_ {{\bf H}} $ converges weakly  in $C([0,T] ^{d} )$  to the $d$-parameter stochastic process $(X_{{\bf t}})_{{\bf t}\geq 0}$ defined by
\begin{equation}
\label{30s-3}
X _{{\bf t}}= \langle {\bf t}\rangle _{d} \frac{1}{\sqrt{q!}}H_{q}(Z) 
\end{equation}
where $Z\sim N(0,1)$ and $H_{q}$ is the $q$th Hermite polynomial (see (\ref{Hermite-poly})).

\item Assume ${\bf H}_{A_{k}}\to \left( \frac{1}{2},..., \frac{1}{2}\right) \in \mathbb{R} ^{k}$. Assume that the parameters $H_{j}, j\in \overline{A}_{k}$ are fixed.  Then the process $ Z^{q,d}_ {{\bf H}}$ converges weakly  in $C([0,T] ^{d} )$ to a $d$-parameter centered Gaussian process $(X ({\bf t}))_{{\bf t}\geq 0}$   with covariance \begin{equation}\label{28s-1}
 \mathbf{E} X_{{\bf t}} X _{{\bf s}} = \left(\prod_{a\in A_{k}} \left( t ^{(a) }\wedge s^{(a)} \right) \right) \left( \prod_{b\in \overline{A} _{k}} R_{H_{b}}(t ^{(b)}, s ^{(b)})\right)
\end{equation}
with $R_{H_{b}}$ defined in (\ref{covh}).

\item Assume ${\bf H}_{A_{k}}\to \left( \frac{1}{2},.., \frac{1}{2} \right) \in \mathbb{R} ^{k} $ and ${\bf H} _{B_{p}}\to
\left( 1,..,1 \right) \in \mathbb{R} ^{p}$.  Assume that  the  $H_{j}$ with $ j\in \{1,2,.., d\} \setminus (A_{k}\cup B_{p}) $ are fixed. 
Then the process $ Z^{q,d }_ {{\bf H}} $ converges weakly  in $C([0,T] ^{d} )$ to a $d$-parameter Gaussian process $(X ({\bf t}))_{{\bf t}\geq 0}$   with covariance
\begin{equation}
\label{28s-4} \mathbf{E} X_{{\bf t}} X _{{\bf s}} =\left( \prod_{ a\in A_{k} }( t ^{(a)} \wedge s^{(a)}) \right) \left( \prod _{b\in B_{p}} t ^{(b)} s ^{(b)} \right) \left( \prod _{c\in \overline{A_{k}\cup B_{p}} }R_{H_{c}} (t^{(c)}, s^{(c)})\right).
\end{equation}

\end{enumerate}
\end{theorem}

We will use the above result in order to get the limit behavior with respect to the Hurst parameter of the Hermite Wiener integral.

\section{Convergence of the Wiener-Hermite integrals with respect to the Hurst parameter}

Let us start the analysis of the behavior of the Wiener-Hermite integral ({\ref{hw}) when the components of the self-similarity index ${\bf H}$ tends to their extreme values.  As mentioned above, we will separate our study into two cases: at least one component of ${\bf H}$ converges to 1 (and no component tends to $\frac{1}{2}$) and at least one component of ${\bf H}$ converges to one-half.

\subsection{Convergence around 1}

We need to introduce new spaces for the deterministic integrand in (\ref{hw}). Working on these spaces will ensure the convergence of the Hermite-Wiener integral.

Let $A_{k}$ be as in (\ref{not}) and assume $1\leq k <d$. We introduce the space $\mathcal{H} _{\overline{A}_{k}} $ of measurable functions $f: \mathbb{R} ^{d} \to \mathbb{R} $ such that 
\begin{eqnarray}
&&{\Vert f \Vert} _{\mathcal{H}_{\overline{A}_{k}}}:= \\
&&\sum_{j=1} ^{k} 
\int_{\mathbb{R} ^{j} }d{\bf u} _{A_{j}} \left| \int_{\mathbb{R} ^{d-j}} d{\bf v} _{\overline{A}_{j}}  \int_{\mathbb{R}^{d-j}}d{\bf w} _{\overline{A}_{j}} \vert f({\bf u}_{A_{j}}, {\bf v}_{\overline{A}_{j}})\vert \cdot \vert f({\bf u}_{A_{j}}, {\bf w}_{\overline{A}_{j}})\vert 
\vert {\bf v} _{\overline{A}_{j}}-{\bf w} _{\overline{A}_{j}}\vert ^ {2{\bf H}_{\overline{A}_{j}}-2}\right| ^{\frac{1}{2}} \nonumber \\
&=&\sum_{j=1} ^{k} 
\int_{\mathbb{R} ^{j} }d{\bf u} _{A_{j}} \Vert f({\bf u}_{A_{j}}, \cdot )\Vert _{ {\mathcal H}_{ {\bf H}_{\overline{A}_{j}}}} <\infty\label{ha}
\end{eqnarray}
with the norm $\Vert \cdot \Vert _{ {\mathcal H}_{ {\bf H}_{\overline{A}_{j}}}}$ defined in (\ref{normahb}). Notice that for $f\in \mathcal{H}_{\overline{A}_{k}}$, the integral 
\begin{equation}\label{27a-4}
\int_{\mathbb{R} ^{k} }d{\bf u} _{A_{k}} \int_{\mathbb{R} ^{d-k}} dZ^{q,d}_{{\bf H}} (u_{\overline{A}_{k}}) f({\bf u})
\end{equation}
is well-defined in $L^{1}(\Omega).$ Indeed,
{
\begin{eqnarray*}
&&\mathbf{E}\left| \int_{\mathbb{R} ^{k} }d{\bf u} _{A_{k}} \int_{\mathbb{R} ^{d-k}} dZ^{q,d}_{{\bf H}} ({\bf u}_{\overline{A}_{k}}) f({\bf u})\right|\\
 &\leq&  \int_{\mathbb{R} ^{k} }d{\bf u} _{A_{k}}  \mathbf{E}\left|\int_{\mathbb{R} ^{d-k}} dZ^{q,d}_{{\bf H}} ({\bf u}_{\overline{A}_{k}}) f({\bf u})\right|
\leq  \int_{\mathbb{R} ^{k} }d{\bf u} _{A_{k}}  \left( \mathbf{E}\left|\int_{\mathbb{R} ^{d-k}} dZ^{q,d}_{{\bf H}} ({\bf u}_{\overline{A}_{k}}) f({\bf u})\right|^{2} \right) ^{\frac{1}{2}} \\
&=& \left( {\bf H}_{\overline{A}_{k}} (2{\bf H}_{\overline{A}_{k}}-{\bf 1}) \right) ^{\frac{1}{2}} \int_{\mathbb{R} ^{k} }d{\bf u} _{A_{k}} \\
&&\left| \int_{\mathbb{R} ^{d-k}} d{\bf v} _{\overline{A}_{k}}  \int_{\mathbb{R}^{d-k}}d{\bf w} _{\overline{A}_{k}} \vert f({\bf u}_{A_{k}}, {\bf v}_{\overline{A}_{k}})\vert \cdot \vert f({\bf u}_{A_{k}}, {\bf w}_{\overline{A}_{k}})\vert 
\vert {\bf v} _{\overline{A}_{j}}-{\bf w} _{\overline{A}_{j}}\vert ^ {2{\bf H}_{\overline{A}_{j}}-2}\right| ^{\frac{1}{2}} \\
& \leq & \left( {\bf H}_{\overline{A}_{k}} (2{\bf H}_{\overline{A}_{k}}-{\bf 1}) \right) ^{\frac{1}{2}}\sum_{j=1} ^{k} 
\int_{\mathbb{R} ^{j} }d{\bf u} _{A_{j}}\\
&& \left| \int_{\mathbb{R} ^{d-j}} d{\bf v} _{\overline{A}_{j}}  \int_{\mathbb{R}^{d-j}}d{\bf w} _{\overline{A}_{j}} \vert f({\bf u}_{A_{j}}, {\bf v}_{\overline{A}_{j}})\vert \cdot \vert f({\bf u}_{A_{j}}, {\bf w}_{\overline{A}_{j}})\vert 
\vert {\bf v} _{\overline{A}_{j}}-{\bf w} _{\overline{A}_{j}}\vert ^ {2{\bf H}_{\overline{A}_{j}}-2}\right| ^{\frac{1}{2}} \\
&= & \left( {\bf H}_{\overline{A}_{k}} (2{\bf H}_{\overline{A}_{k}}-{\bf 1}) \right) ^{\frac{1}{2}} \Vert f\Vert _{\mathcal{H}_{\overline{A}_{k}}}<\infty.
\end{eqnarray*}
}
If $k=d$, we define  
$\mathcal{H}_{\overline{A}_{k}}=\mathcal{H}_{\overline{A}_{d}}$ to be the set of measurable functions $f: \mathbb{R} ^{d} \to \mathbb{R} $ such that 
\begin{eqnarray}
&&\Vert f \Vert  _{\mathcal{H}_{\overline{A}_{k}}}:=\Vert f\Vert  _{ L ^{1}(\mathbb{R}^{d})}\nonumber\\
&&+ \sum_{j=1} ^{d-1} 
\int_{\mathbb{R} ^{j} }d{\bf u} _{A_{j}} \left| \int_{\mathbb{R} ^{d-j}} d{\bf v} _{\overline{A}_{j}}  \int_{\mathbb{R}^{d-j}}d{\bf w} _{\overline{A}_{j}} \vert f({\bf u}_{A_{j}}, {\bf v}_{\overline{A}_{j}})\vert \cdot \vert f({\bf u}_{A_{j}}, {\bf w}_{\overline{A}_{j}})\vert 
\vert {\bf v} _{\overline{A}_{j}}-{\bf w} _{\overline{A}_{j}}\vert ^ {2{\bf H}_{\overline{A}_{j}}-2}\right| ^{\frac{1}{2}}  \nonumber\\
&:=&\Vert f\Vert  _{ L ^{1}(\mathbb{R}^{d})} +  {\Vert f \Vert} _{\mathcal{H}_{\overline{A}_{d-1}}} <\infty.\label{ha1}
\end{eqnarray}

\begin{remark}
Notice that the order of integration in (\ref{27a-4}) is important. That is, the integral
$$ \int_{\mathbb{R} ^{d-k}} dZ^{q,d}_{{\bf H}} (u_{\overline{A}_{k}}) \int_{\mathbb{R} ^{k} }d{\bf u} _{A_{k}} f({\bf u}) $$
is not necesarily well-defined for $f\in \mathcal{H}_{\overline{A}_{k}}.$
\end{remark}

We have the following non-central limit theorem.
\begin{prop}\label{pp1} 
Let $A_{k}$ be as in (\ref{not})  and assume  $f\in \mathcal{H}_{\overline{A}_{k}}\cap  \vert \mathcal{H} _{{\bf H}}\vert$.

\begin{itemize}

\item Assume $1\leq k<d$  and

$$ {\bf H}_{A_{k}} \to (1,..,1) \in \mathbb{R} ^{k}  \mbox{ and } {\bf H}_{\overline{A}_{k}}\in \left(\frac{1}{2}, 1\right) ^{d-k} \mbox{ is fixed. } $$
 Then the family of random variables $\left( X^{{\bf H}}, {\bf H} \in \left( \frac{1}{2}, 1\right) ^{d} \right) $
\begin{equation}\label{xh}
X^{{\bf H}}:=\int_{\mathbb{R}^{d}}f({\bf u}) d Z^{q,d}_{{\bf H}}({\bf u})
\end{equation}
converges in distribution to the random variable
\begin{equation}
\label{x}X:= \int_{\mathbb{R} ^{d}} f ( u^{(1)},.., u^{(d) }) dZ ^{q,d-k} _{\overline {A}_{k}} ({\bf u}_{\overline{A}_{k}})d{\bf u} _{A_{k}} = \int_{\mathbb{R} ^{k}}\left(    \int_{\mathbb{R} ^{d-k}}  f({\bf u}_{A_{k}}, {\bf u} _{\overline{A}_{A_{k}}})dZ ^{q,d-k} _{\overline {A}_{k}} ({\bf u}_{\overline{A}_{k}})\right) d{\bf u} _{A_{k}} .
\end{equation}

\item Assume  $k=d$ and
$${\bf H}\to (1,..,1) \in \mathbb{R} ^{d}.$$
 Then the limit in distribution of the  family  $\left( X^{{\bf H}}, {\bf H} \in \left( \frac{1}{2}, 1\right) ^{d} \right) $ given by  (\ref{xh}) is
$$\int_{\mathbb{R} ^{d}} f ( u^{(1)},.., u^{(d) }) d{\bf u} \frac{1}{\sqrt{q!}}H_{q} (Z)$$
with $Z\sim N(0,1)$ and $H_{q}$ the Hermite polinomial of degree $q$ (\ref{Hermite-poly}).

\end{itemize}
\end{prop}
{\bf Proof: } We will check the convergence of the characteristic function of $X^{{\bf H}}$. That is, we will show that for every $\alpha \in \mathbb{R} $,

$$\mathbf{E} e ^{i \alpha X ^{{\bf H}} }\to _{{\bf H}_{A_{k}} \to (1,..,1)\in \mathbb{R} ^{k} } \mathbf{E} e ^{i\alpha X}.$$

The idea is to approximate first  $X ^{{\bf }}$ by a sequence of random variables that can be written in terms of the linear combinaisons of $ Z ^{q, d}_{{\bf H}}$ and to use the result in Theorem \ref{tat}. Consider a sequence of step functions 
$$f_{n}({\bf u}) = \sum_{l=1} ^{n}  a_{l} 1_{ ( {\bf t}_{l}, {\bf t}_{l+1}]} ({\bf u}) = \sum_{l=1} ^{n}  a_{l} 1_{ ( t_{l}^{(1)}, t_{l+1} ^{(1)}]}(u^{(1)} )....1_{ ( t_{l}^{(d)}, t_{l+1} ^{(d)}]}(u^{(d)} )$$
(where we used again the notation  ${\bf u}= (u ^{(1)},.., u^{(d)}) $ and $ {\bf t}_{l} =(t_{l} ^{(1)},.., t_{l}^{(d)})$ for $l=1,..,n$) such that

\begin{equation}
\label{24i-1}
 \Vert f_{n} -f\Vert _{ \mathcal{H}_{\overline{A}_{k}}} \to_{n\to \infty}0\mbox{ and } \Vert f_{n} -f\Vert _{ \vert \mathcal{H}_{{\bf H}}\vert}\to_{n\to \infty}0.
\end{equation}
The choice  of such a sequence $ (f_{n})_{n\geq 1}$ is possible because for any positive function $f\in  \mathcal{H}_{\overline{A}_{k}}\cap  \vert \mathcal{H} _{{\bf H}}\vert$, there exists an increasing  sequence of step functions in ${f_{n}}\in \mathcal{H}_{\overline{A}_{k}}\cap  \vert \mathcal{H} _{{\bf H}}\vert$ which converges poinwise to $f$ and satisfies $\vert f_{n}-f\vert \leq \vert f\vert$, and by dominated convergence theorem, it converges in $\mathcal{H}_{\overline{A}_{k}}$ and in $ \vert \mathcal{H} _{{\bf H}}\vert$. Then, we use the fact that a general function can be decomposed into its positive and negative parts.

Consider the Hermite Wiener integral of $f_{n}$ with respect to the Hermite sheet 
\begin{equation*}
X ^{n, {\bf H}}= \int_{\mathbb{R} ^{d} } f_{n}({\bf u}) dZ ^{q,d}_{{\bf H}} ({\bf u})=\sum _{j=1} ^{n} a_{l} (\Delta Z ^{q,d} _{{\bf H}} )( ({\bf t}_{l}, {\bf t}_{l+1} ])
\end{equation*}
with $\Delta Z ^{q,d} _{{\bf H}} $ given by (\ref{mi}). Then we know  from \cite{DCT}, Section 3 that $X ^{n, {\bf H}}$ converges in $L ^{2} (\Omega)$ to $X ^{ {\bf H}}$  if $f_{n}$ converges to $f$ in $\vert \mathcal{H}_{{\bf H}}\vert$ due to  the isometry of the Hermite Wiener integral (\ref{isoHW}). So we have 
$$X ^{n, {\bf H}} \to _{n\to \infty} X^{{\bf H}}:= \int_{\mathbb{R} ^{d}} f({\bf s}) dZ ^{q,d}_{{\bf H}} ({\bf s}) \mbox{ in } L ^{2}(\Omega).$$

Consequently, we can write
\begin{equation}
\label{27a-2}
\lim _{{\bf H}_{A_{k}} \to (1,..,1)\in \mathbb{R} ^{k} }\mathbf{E} e ^{i \alpha X ^{{\bf H}} } = \lim _{{\bf H}_{A_{k}} \to (1,..,1)\in \mathbb{R} ^{k} }\lim _{n\to \infty} \mathbf{E} e ^{i \alpha X ^{n, {\bf H}} }.
\end{equation}

Now, we aim at exchanging the two limits above. Recall that if $f_{j}, j\geq 1$ is a sequence of functions on $D\subset \mathbb{R}$ converging uniformly to $f$ on $D$ and if $ a$ is a limit point for $D$, then $ \lim _{j\to \infty} \lim _{x\to a} f_{j} (x)= \lim _{x\to a} f(x) $ provided that $\lim _{x\to a} f(x), \lim _{x\to a} f_{j}(x)$ exist. Therefore it suffices to show that $\mathbf{E} e ^{i \alpha X ^{n, H} }$ converges uniformly with respect to ${\bf H }_{A_{k}} $ to $\mathbf{E} e ^{i \alpha X ^{H} }.$

By the mean value theorem

$$ \left| \mathbf{E} e ^{i \alpha X ^{n, {\bf H}} }-\mathbf{E} e ^{i \alpha X ^{ {\bf H}} }\right| \leq \vert \alpha \vert \mathbf{E} \left| X ^{n,{\bf H}} -X ^{{\bf H}}\right| \leq \vert \alpha \vert \left( \mathbf{E} \left| X ^{n,{\bf H}} -X ^{{\bf H}}\right|^{2} \right) ^{\frac{1}{2}}.  $$

Thus, in order to invert the limits in (\ref{27a-2}),  it suffices to show that for some $\eps >0$

$$\sup _{{\bf H} _{A_{k}} \in [\frac{1}{2}+ \eps, 1]^{k} }\mathbf{E} \left| X ^{n,{\bf H}} -X ^{{\bf H}}\right|^{2} \to _{n\to \infty } 0$$
that is proved in Lemma \ref{ll1} below. The relation (\ref{27a-2}) becomes

\begin{equation}
\label{20s-1}
\lim _{{\bf H}_{A_{k}} \to (1,..,1)\in \mathbb{R} ^{k} }\mathbf{E} e ^{i \alpha X ^{{\bf H}} } =\lim _{n\to \infty}\lim _{{\bf H}_{A_{k}} \to (1,..,1)\in \mathbb{R} ^{k} }  \mathbf{E} e ^{i \alpha X ^{n, {\bf H}} }.
\end{equation}

Assume $k<d$. Since, from Theorem \ref{tat}
$Z ^{q, d}_{{\bf H}} $ converges weakly to the process $\left( U_{{\bf t}}\right) _{{\bf t}\geq 0}$ given by 
$$U_{{\bf t}} =\langle {\bf t} \rangle _{A_{k}} Z ^{q, d-k} _{{ \bf H} _{ \overline{A}_{k}}}({\bf t} _{\overline{A}_{k}})$$ 
it follows from (\ref{20s-1}) that
\begin{eqnarray}
\lim _{{\bf H}_{A_{k}} \to (1,..,1)\in \mathbb{R} ^{k} }\mathbf{E} e ^{i \alpha X ^{{\bf H}} }&=& \lim_{n\to \infty}  \lim _{{\bf H}_{A_{k}} \to (1,..,1)\in \mathbb{R} ^{k} }\mathbf{E} e ^{i\alpha  \sum _{l=1} ^{n} a_{l} (\Delta Z ^{q,d} _{{\bf H}} )( ({\bf t}_{l}, {\bf t}_{l+1} ])}\nonumber \\
&=& \lim_{n\to \infty}\mathbf{E} e ^{i \alpha \sum_{l=1}^{n}  a_{l} (\Delta  U) ( ({\bf t}_{l}, {\bf t}_{l+1} ])}.\label{6o-1}
\end{eqnarray}

At this point we need to  study the convergence as $n\to \infty$ of the sequence
\begin{equation}
\label{20s-2}
X^{n}:=\sum_{l=1}^{n}  a_{l} (\Delta  U) ( ({\bf t}_{l}, {\bf t}_{l+1} ])
\end{equation}
as $n\to \infty$. If $A_{k}=\{j_{1},.., j_{k}\}$, let us use the notation

$$({\bf t}_{l}, {\bf t}_{l+1} ]_{A_{k}}= (t_{l} ^{(j_{1})}, t_{l+1} ^{(j_{1})}]\times....\times (t_{l} ^{(j_{k})}, t_{l+1} ^{(j_{k})}].$$
Then it is not difficult to see that 
$$(\Delta  U) ( ({\bf t}_{l}, {\bf t}_{l+1} ])= (\Delta \langle {\bf t}\rangle _{A_{k}}) ({\bf t}_{l}, {\bf t}_{l+1} ]_{A_{k}} (\Delta Z ^{q,d-k} _{\overline {A}_{k}}) ({\bf t}_{l}, {\bf t}_{l+1} ]_{\overline{A}_{k}}.$$
and therefore the sequence (\ref{20s-2}) can be expressed as follows
\begin{eqnarray*}
X^{n}&=&\sum_{l=1}^{n}  a_{l} (\Delta  U) ( ({\bf t}_{l}, {\bf t}_{l+1} ]]=\sum_{l=1}^{n}  a_{l} (\Delta \langle {\bf t}\rangle _{A_{k}}) ({\bf t}_{l}, {\bf t}_{l+1} ]_{A_{k}} (\Delta Z ^{q,d-k} _{{\bf H}_{\overline {A}_{k}}}) ({\bf t}_{l}, {\bf t}_{l+1} ]_{\overline{A}_{k}}\\
&=& \int_{\mathbb{R} ^{d}} f_{n} ( u^{(1)},.., u^{(d) }) d{\bf u} _{A_{k}} d Z ^{q,d-k} _{{\bf H}_{\overline {A}_{k}}} ({\bf u}_{\overline{A}_{k}}). 
\end{eqnarray*}

Now, we show that 

\begin{equation}\label{6o-2}
X ^{n} \to _{n \to \infty} X \mbox{ in } L ^{1} (\Omega)
\end{equation}
where the random variable $X$ is given by (\ref{x}). We have

\begin{eqnarray*}
\mathbf{E} \vert X^{n}-X\vert &=&\mathbf{E} \left|   \int_{\mathbb{R} ^{k}} d{\bf u}_{A_{k}} \int_{\mathbb{R}^{d-k}} dZ^{q,d}_{{\bf H}} ({\bf u} _{\overline{A}_{k}}  ) (f_{n}( {\bf u}_{A_{k}},  {\bf u} _{\overline{A}_{k}})-f( {\bf u}_{A_{k}},  {\bf u} _{\overline{A}_{k}})) \right|\\
&\leq &  \int_{\mathbb{R} ^{k}} d{\bf u}_{A_{k}} \mathbf{E}\left|  \int_{\mathbb{R}^{d-k}} dZ^{q,d}_{{\bf H}} ({\bf u} _{\overline{A}_{k}}  ) (f_{n}( {\bf u}_{A_{k}},  {\bf u} _{\overline{A}_{k}})-f( {\bf u}_{A_{k}},  {\bf u} _{\overline{A}_{k}}))  \right| \\
&\leq &\int_{\mathbb{R} ^{k}} d{\bf u}_{A_{k}} \left( \mathbf{E}\left|  \int_{\mathbb{R}^{d-k}} dZ^{q,d}_{{\bf H}} ({\bf u} _{\overline{A}_{k}}  ) (f_{n}( {\bf u}_{A_{k}},  {\bf u} _{\overline{A}_{k}})-f( {\bf u}_{A_{k}},  {\bf u} _{\overline{A}_{k}}))  \right| ^{2} \right) ^{\frac{1}{2}}\\
&{=}& {\left( {\bf H}_{\overline{A}_{k}} (2{\bf H}_{\overline{A}_{k}}-{\bf 1}) \right)^{\frac{1}{2}}} \int_{\mathbb{R} ^{k}} d{\bf u}_{A_{k}}\left|  \int_{\mathbb{R}^{d-k}} \int_{\mathbb{R}^{d-k}} d{\bf v}_{\overline{A}_{k}}d{\bf w}_{\overline{A}_{k}} \vert {\bf v}_{\overline{A}_{k}}-{\bf w}_{\overline{A}_{k}}\vert ^{2{\bf H}_{\overline{A}_{k}}-2}\right.\\
&&\left. \times   \left( f_{n}( {\bf u}_{A_{k}},  {\bf v} _{\overline{A}_{k}})-f( {\bf u}_{A_{k}},  {\bf v} _{\overline{A}_{k}})\right) \left( f_{n}( {\bf u}_{A_{k}},  {\bf w} _{\overline{A}_{k}})-f( {\bf u}_{A_{k}},  {\bf w} _{\overline{A}_{k}})\right) \right| ^{\frac{1}{2}}\\
&\leq & {\left( {\bf H}_{\overline{A}_{k}} (2{\bf H}_{\overline{A}_{k}}-{\bf 1}) \right)^{\frac{1}{2}}} \Vert f_{n}-f\Vert _{{\mathcal{H}_{\overline{A}_{k}}}}\to_{n\to \infty} 0
\end{eqnarray*}
where the last convergence comes from (\ref{24i-1}). We obtain from (\ref{6o-1}) and (\ref{6o-2})
$$\lim _{{\bf H}_{A_{k}} \to (1,..,1)\in \mathbb{R} ^{k} }\mathbf{E} e ^{i \alpha X ^{{\bf H}} }= \lim _{n \to \infty} \mathbf{E} e ^{i\alpha X ^{n}} = \mathbf{E} e ^{i \alpha X}$$
and the proof is complete for $1\leq k<d$.

If $k=d$, the proof is  similar. We know that the process
$Z ^{q, d}_{{\bf H}} $ converges weakly in $C[0,T]$ to the process 
$$\langle {\bf t}\rangle_{d} \frac{1}{\sqrt{q!}} H_{q}(Z).$$
Using the same lines as above, we get
$$\lim _{{\bf H}\to (1,..,1)\in \mathbb{R} ^{d} }\mathbf{E} e ^{i \alpha X ^{{\bf H}} }=\lim _{n \to \infty} \mathbf{E} e ^{i\alpha X ^{n}} $$
and in  this case the sequence (\ref{20s-2}) becomes
$$X^{n}=  \sum_{i=1} ^{n} (\Delta \langle {\bf t}\rangle_{d} ) [{\bf t}_{l}, {\bf t}_{l+1}]  \frac{1}{\sqrt{q!}} H_{q}(Z) = \int_{\mathbb{R}} f_{n}({\bf u}) d{\bf u} \frac{1}{\sqrt{q!}} H_{q}(Z)$$
{
Clearly, by (\ref{24i-1})
$$\mathbf{E}\vert X^{n}- \int_{{\mathbb{R}^{d}}} f({\bf u}) d{\bf u} \frac{1}{\sqrt{q!}} H_{q}(Z)\vert  \leq \left( \int_{\mathbb{\mathbb{R}^{d}}}\vert f_{n}({\bf u}) - f({\bf u}) \vert d{\bf u}\right) \frac{1}{\sqrt{q!}} H_{q}(Z)\to _{n\to \infty}0.$$
using the definition of the norm in  $\mathcal{H}_{\overline{A}_{k}}$ for $k=d$.
Then
$$\lim _{n\to \infty} \mathbf{E} e ^{i \alpha X ^{n} }=\mathbf{E} e ^{i\alpha (\int_{\mathbb{R} ^{d}} f({\bf u}) d{\bf u})\frac{1}{\sqrt{q!}} H_{q}(Z) }.$$
\qed 

}

The below lemma has been needed in the proof of Proposition \ref{pp1}.
\begin{lemma}\label{ll1} Let $A_{k}$ be as in (\ref{not}) with $1\leq k\leq d$. Assume $f\in \mathcal{H}_{\overline{A}_{k}}\cap  \vert \mathcal{H} _{{\bf H}}\vert$ and consider a sequence $ (f_{n})_{n\geq 1} $ of step functions on $\mathbb{R} ^{d}$ such that (\ref{24i-1}) holds true.
Let
$$X^{n,{\bf H}}= \sum _{l=1} ^{n} a_{l} (\Delta Z ^{q,d} _{{\bf H}} )( ({\bf t}_{l}, {\bf t}_{l+1} ]).$$
Then for every $\eps>0$ small enough
$$\sup _{{\bf H} _{A_{k}} \in [\frac{1}{2}+ \eps, 1]^{k} }\mathbf{E} \left| X ^{n,{\bf H}} -X ^{{\bf H}}\right|^{2} \to _{n\to \infty } 0.$$
\end{lemma}
{\bf Proof: }
From the isometry property (\ref{isoHW}) and from (\ref{24i-1}) we have for every ${\bf H} \in (\frac{1}{2}, 1) ^{d} $, 
\begin{equation}\label{24i-4}
\mathbf{E}\left| X ^{n,{\bf H}} -X ^{{\bf H}}\right|^{2}\to 0.
\end{equation}
Let us show that the above convergence is uniform with respect to ${\bf H} _{A_{k}} \in [\frac{1}{2}+ \eps, 1]^{k} $. By  (\ref{isoHW}),
\begin{eqnarray}
\mathbf{E} \left| X ^{n,{\bf H}} -X ^{{\bf H}}\right|^{2}&=& {\bf H}(2{\bf H}-{\bf 1}) \int_{\mathbb{R} ^{d}} \int_{\mathbb{R} ^{d}}f_{n}({\bf u}) f_{n}({\bf v}) \vert {\bf u}-{\bf v}\vert ^{2{\bf H}-2}d{\bf u} d{\bf v} \nonumber \\
&&-2{\bf H}(2{\bf H}-{\bf 1}) \int_{\mathbb{R} ^{d}} \int_{\mathbb{R} ^{d}}f_{n}({\bf u}) f({\bf v}) \vert {\bf u}-{\bf v}\vert ^{2{\bf H}-2}d{\bf u} d{\bf v}\nonumber \\
&&+ {\bf H}(2{\bf H}-{\bf 1}) \int_{\mathbb{R} ^{d}} \int_{\mathbb{R} ^{d}}f({\bf u}) f({\bf v}) \vert {\bf u}-{\bf v}\vert ^{2{\bf H}-2}d{\bf u} d{\bf v}\nonumber \\
:=G({\bf H}_{A_{k}})\label{20s-4}
 \end{eqnarray}
with the function $G$ considered on the interval $ [\frac{1}{2}+ \eps, 1]^{k}.$ Assume $k<d$. Let  ${\bf 1}(A_{k})=(1,..,1)\in \mathbb{R} ^{k}$. Then  from  (\ref{20s-4}) 
\begin{eqnarray*}
&&G({\bf 1}(A_{k}))\\
&=&{\bf H }_{\overline{A}_{k}} (2{\bf H }_{\overline{A}_{k}} -{\bf 1})\int_{\mathbb{R}^{k}} d{\bf u}_{{A_{k}}}\int_{\mathbb{R}^{k}} d{\bf v}_{{A_{k}}}\int_{\mathbb{R}^{d-k}} d{\bf u}_{{\overline{A}_{k}}}
\int_{\mathbb{R}^{d-k}} d{\bf v}_{{\overline{A}_{k}}}f_{n} ( {\bf u}_{{A_{k}}}, {\bf u}_{{\overline{A}_{k}}})
f_{n} ( {\bf v}_{{A_{k}}}, {\bf v}_{{\overline{A}_{k}}}){\vert {\bf u}_{{\overline{A}_{k}}} - {\bf v}_{{\overline{A}_{k}}} \vert^{2{\bf H }_{\overline{A}_{k}} -2}} \nonumber\\
&&-2{\bf H }_{\overline{A}_{k}} (2{\bf H }_{\overline{A}_{k}} -{\bf 1})\int_{\mathbb{R}^{k}} d{\bf u}_{{A_{k}}}\int_{\mathbb{R}^{k}} d{\bf v}_{{A_{k}}}\int_{\mathbb{R}^{d-k}} d{\bf u}_{{\overline{A}_{k}}}
\int_{\mathbb{R}^{d-k}} d{\bf v}_{{\overline{A}_{k}}}f_{n} ( {\bf u}_{{A_{k}}}, {\bf u}_{{\overline{A}_{k}}})
f( {\bf v}_{H_{A_{k}}}, {\bf v}_{{\overline{A}_{k}}}){\vert {\bf u}_{{\overline{A}_{k}}} - {\bf v}_{{\overline{A}_{k}}} \vert^{2{\bf H }_{\overline{A}_{k}} -2}}\nonumber\\
&&+{\bf H }_{\overline{A}_{k}} (2{\bf H }_{\overline{A}_{k}} -{\bf 1})\int_{\mathbb{R}^{k}} d{\bf u}_{{A_{k}}}\int_{\mathbb{R}^{k}} d{\bf v}_{{A_{k}}}\int_{\mathbb{R}^{d-k}} d{\bf u}_{{\overline{A}_{k}}}
\int_{\mathbb{R}^{d-k}} d{\bf v}_{{\overline{A}_{k}}}f ( {\bf u}_{{A_{k}}}, {\bf u}_{{\overline{A}_{k}}})
f ( {\bf v}_{{A_{k}}}, {\bf v}_{{\overline{A}_{k}}}){\vert {\bf u}_{{\overline{A}_{k}}} - {\bf v}_{{\overline{A}_{k}}} \vert^{2{\bf H }_{\overline{A}_{k}} -2}}\nonumber
\end{eqnarray*}
and this can be written
\begin{eqnarray}
&&G({\bf 1}(A_{k}))\nonumber \\
&=&\int_{\mathbb{R}^{k}} d{\bf u}_{{A_{k}}}\int_{\mathbb{R}^{k}} d{\bf v}_{{A_{k}}}\langle(f_{n}- f)({\bf u}_{{A_{k}}}, \cdot), (f_{n}-f)({\bf v}_{{A_{k}}}, \cdot)\rangle _{ \mathcal{H} _{H_{\overline{A}_{k}}}}\nonumber \\
&\leq &\int_{\mathbb{R}^{k}} d{\bf u}_{{A_{k}}}\int_{\mathbb{R}^{k}} d{\bf v}_{{A_{k}}}\Vert (f_{n}- f)({\bf u}_{{A_{k}}}, \cdot)\Vert  _{ \mathcal{H} _{H_{\overline{A}_{k}}}}\Vert (f_{n}- f)({\bf v}_{{A_{k}}}, \cdot)\Vert  _{ \mathcal{H} _{H_{\overline{A}_{k}}}}\nonumber \\
&=& \left(\int_{\mathbb{R}^{k}} d{\bf u}_{{A_{k}}}\Vert (f_{n}- f)({\bf u}_{{A_{k}}}, \cdot)\Vert  _{ \mathcal{H} _{H_{\overline{A}_{k}}}}\right) ^{2}\nonumber \\
&=&{{\bf H }_{\overline{A}_{k}} (2{\bf H }_{\overline{A}_{k}} -{\bf 1}) \left[\int_{\mathbb{R} ^{k} }d{\bf u} _{A_{k}} \left| \int_{\mathbb{R} ^{d-k}} d{\bf v} _{\overline{A}_{k}}  \int_{\mathbb{R}^{d-k}}d{\bf w} _{\overline{A}_{k}} \vert f({\bf u}_{A_{k}}, {\bf v}_{\overline{A}_{k}})\vert \cdot \vert f({\bf u}_{A_{k}}, {\bf w}_{\overline{A}_{k}})\vert 
\vert {\bf v} _{\overline{A}_{k}}-{\bf w} _{\overline{A}_{k}}\vert ^ {2{\bf H}_{\overline{A}_{k}}-2}\right| ^{\frac{1}{2}} \right]^{2}}\nonumber\\
&\leq &{\bf H }_{\overline{A}_{k}} (2{\bf H }_{\overline{A}_{k}} -{\bf 1}) \Vert f_{n}-f\Vert{^{2}} _{\mathcal{H}_{\overline{A}_{k}}}
\label{24i-3}
\end{eqnarray}
where we used the  definition (\ref{ha}).

Now, the function $G$ is continuous on $ [\frac{1}{2}+ \eps, 1]^{k}$ so there exists ${\bf H_{0} } = ( H_{0,1},.., H_{0, k})\in  [\frac{1}{2}+ \eps, 1]^{k}$ such that 

$$
\sup _{{\bf H} _{A_{k}} \in [\frac{1}{2}+ \eps, 1]^{k} }{G({\bf H} _{A_{k}})= G({\bf H}_{0})}.$$

If ${\bf H}_{0}= {\bf 1}(A_{k})$, then the conclusion follows from (\ref{24i-3}) and the assumption (\ref{24i-1}).  If ${\bf H}_{0} $ has the form 
$${\bf H}_{0} =(1,.., 1, H_{0, j+1},..., H_{0, k})$$
with $j<k$ then a similar calculation to (\ref{24i-3}) shows that 
\begin{eqnarray}
&&G({\bf H}_{0} )\nonumber\\
 &\leq &  {{\bf H }_{\overline{A}_{j}} (2{\bf H }_{\overline{A}_{j}} -{\bf 1})\left[\int_{\mathbb{R} ^{j} }d{\bf u} _{A_{j}} \left| \int_{\mathbb{R} ^{d-j}} d{\bf v} _{\overline{A}_{j}}  \int_{\mathbb{R}^{d-j}}d{\bf w} _{\overline{A}_{j}} \vert f({\bf u}_{A_{j}}, {\bf v}_{\overline{A}_{j}})\vert \cdot \vert f({\bf u}_{A_{j}}, {\bf w}_{\overline{A}_{j}})\vert 
\vert {\bf v} _{\overline{A}_{j}}-{\bf w} _{\overline{A}_{j}}\vert ^ {2{\bf H}_{\overline{A}_{j}}-2}\right| ^{\frac{1}{2}} \right]^{2}}\nonumber\\
&\leq &{\bf H }_{\overline{A}_{j}} (2{\bf H }_{\overline{A}_{j}} -{\bf 1})\Vert f_{n}-f\Vert{^{2}}_{\mathcal{H}_{\overline{A}_{k}}}\label{6o-3}
\end{eqnarray}
and again $G({\bf H}_{0} ) \to 0$ as ${\bf H}_{A_{k}}\to (1,..,1) \in \mathbb{R} ^{k}$ from (\ref{24i-1}).

Otherwise, if all $H_{0,i}, i=1,..,k$  are in $[\frac{1}{2}+\eps, 1)$, then the conclusion follows from (\ref{24i-4}). 

If $k=d$, the conclusion follows in the same way.
Let $G$ be given by (\ref{20s-4}) and  let ${\bf H_{0} } = ( H_{0,1},.., H_{0, k})\in  [\frac{1}{2}+ \eps, 1]^{d}$ such that 

$$
\sup _{{\bf H} \in [\frac{1}{2}+ \eps, 1]^{g} }{G({\bf H})= G({\bf H}_{0})}.$$

 If ${\bf H}_{0}=(1,.., 1) \in \mathbb{R} ^{d}$, notice that in this case ${G({\bf 1}_{d})}=G(1,..., 1)= \Vert f_{n}-f\Vert ^{2} _{L ^{1}(\mathbb{R} ^{d})}\to _{n\to \infty} 0$.  If ${\bf H}_{0} $ has the form 
$${\bf H}_{0} =(1,.., 1, H_{0, j+1},..., H_{0, d})$$
with $j<d$ then $G({\bf H}_{0})$ satisfies (\ref{6o-3}) and consequently it converges to zero from  the assumption (\ref{24i-1}). Il all  components of ${\bf H}_{0}$ are strictly contained in the interval $\left(\frac{1}{2}, 1\right)$, then we conclude by (\ref{24i-4}). 
\qed

\subsection{ Convergence around $\frac{1}{2}$}

In this section, we will study the convergence in distribution of the Hermite Wiener integral (\ref{xh}) when at least one Hurst index converges to one half.  Actually, we will  assume (recall notation (\ref{not}) from the previous section)
$$ {\bf H}_{A_{k}} \to \left(\frac{1}{2},..., \frac{1}{2}\right) \in \mathbb{R} ^{k} $$
and 
$${\bf H}_{B_{p}} \to (1,.., 1) \in \mathbb{R} ^{p}$$
with $1\leq k\leq d, 0\leq p\leq d$  and {$ p+k\leq d. $} Note that $k\geq 1$ means that at least one Hurst parameter converges to $\frac{1}{2}$ while $p\geq 0$ means that some Hurst parameters (possibly zero) converges to 1. 

\vskip0.2cm
We have the following result.

\begin{prop}\label{pp2} Assume $A_{k}$ is as in (\ref{not}) 
and $B_{p}= \{l_{1},.., l_{p} \} \subset \{1,..,d\} $ with $ 0\leq p\leq d, 1\leq k\leq d,  p+k\leq d$ and $A_{k}\cap B_{p} =\emptyset$ (if $p=0$ then $B_{p}=\emptyset. $). 
Let $f\in \vert \mathcal{H}_{{\bf H}}\vert$.  Assume that the following limit exists
\begin{equation}\label{c1-12}
\lim _{ {\bf H}_{A_{k}} \to (\frac{1}{2},..., \frac{1}{2}) \in \mathbb{R} ^{k} } {\bf H} (2{\bf H}-{\bf 1}) \int_{\mathbb{R} ^{d}}\int_{\mathbb{R} ^{d}}f({\bf u}) f({\bf v}) \vert {\bf u}-{\bf v} \vert ^{2{\bf H}-2}d{\bf u} d{\bf v} := \sigma _{f, {\bf H}_{\overline{A}_{k}}} ^{2} 
\end{equation}
and that

\begin{eqnarray}
&&\sup _{{\bf H} _{A_{k}}\in [\frac{1}{2}, 1] ^{k}}  \int_{\mathbb{R} ^{d} } \int_{\mathbb{R} ^{d} } \int_{\mathbb{R} ^{d} } \int_{\mathbb{R} ^{d} } d\mathbf{u}d\mathbf{v}d\mathbf{u}'d\mathbf{v}' f({\bf u}) f({\bf u} ') f({\bf v}) f({\bf v}') \nonumber \\
&&\times \vert {\bf u}- {\bf v}\vert ^{ \frac{2(\bf H -{\bf 1})r}{q}}\vert {\bf u}'- {\bf v}'\vert ^{ \frac{2(\bf{H} -{\bf 1})r}{q}}\vert {\bf u}- {\bf u}'\vert ^{ \frac{2(\bf{H} -{\bf 1})(q-r)}{q}}
\vert {\bf v}- {\bf v}'\vert ^{ \frac{2(\bf{H}-{\bf 1})(q-r)}{q}} <\infty.\label{c2-12}
\end{eqnarray}
If 
$$ {\bf H}_{A_{k}} \to \left(\frac{1}{2},..., \frac{1}{2}\right) \in \mathbb{R} ^{k}, 
  {\bf H}_{B_{p}} \to (1,.., 1) \in \mathbb{R} ^{p} \mbox{ and } {\bf H}_{ \overline{A}_{k} \cup \overline{B}_{p}}\in \left(\frac{1}{2}, 1\right) ^{d-k-p} \mbox{ is fixed}$$
then the Hermite Wiener integral $\int_{\mathbb{R} ^{d} }f({\bf u}) dZ ^{q,d}_{{\bf H}}({\bf u}) $ converges in distribution  to the Gaussian law $N(0, \sigma _{f, {\bf H}_{\overline{A}_{k}}} ^{2})$.

\end{prop}
{\bf Proof: } Recall that by (\ref{hw}), $\int_{\mathbb{R} ^{d}} f({\bf u}) d Z^{q,d}_{{\bf H}} ({\bf u})= I_{q}(Jf) $ with the operator $J$ defined in (\ref{jf}). We can apply the Fourth Moment Theorem to study the normal convergence of (\ref{xh}).

First notice that by assumption  (\ref{c1-12}),  we have 
$$\mathbf{E} \left( \int_{\mathbb{R} ^{d}} f({\bf u}) d Z^{q,d}_{{\bf H}} ({\bf u})\right) ^{2} =  {\bf H} (2{\bf H}-{\bf 1}) \int_{\mathbb{R} ^{d}}\int_{\mathbb{R} ^{d}}f({\bf u}) f({\bf v}) \vert {\bf u}-{\bf v} \vert ^{2{\bf H}-2} d{\bf u} d{\bf v} 
$$
converges to $\sigma _{f, {\bf H}_{\overline{A}_{k}}} ^{2}$.   Therefore, in order to apply the Fourth Moment Theorem (see Theorem \ref{fmt} in the Appendix),  it suffices to show that

$$\Vert Jf \otimes _{r} Jf \Vert _{ L ^{2} ( \mathbb{R} ^{d (2q-2r)})} \to 0 $$
for every ${r=1,..., q-1}$. 

Now, as in the proof of Theorem 3 in \cite{AT} (based on relation (13) in this reference)
\begin{eqnarray*}
( Jf \otimes _{r} Jf)({\bf y}_{1},..,{\bf y}_{2q-2r} ) &=&{\int_{(\mathbb{R} ^{d} ) ^{r }} }Jf ({\bf u}_{1},.., {\bf u}_{r}, {\bf y}_{1},..,{\bf y}_{q-r} )Jf ({\bf u}_{1},.., {\bf u}_{r}, {\bf y}_{q-r-1},..,{\bf y}_{2q-2r} )d{\bf u}_{1}...d{\bf u}_{r} \\
&=&c(\mathbf{H},q)  ^{2} \int_{(\mathbb{R} ^{d} ) ^{r }}d{\bf u}_{1}...d{\bf u} _{r}\\
&& \int_{\mathbb{R} ^{d} } \quad f({\bf u}) \left( \prod _{j=1}^{q-r} (\mathbf{u}-\mathbf{y}_{j})_{+} ^{-\mathbf{\left( \frac{1}{2} + \frac{1-H}{q} \right)} } \right) \left( \prod _{j=1}^{r} (\mathbf{u}-\mathbf{u}_{j})_{+} ^{-\mathbf{\left( \frac{1}{2} + \frac{1-H}{q} \right)} } \right) d\mathbf{u}\\
&&\times \int_{\mathbb{R} ^{d} } \quad f({\bf v})\left( \prod _{j=q-r+1}^{2q-2r} (\mathbf{v}-{\mathbf{y}_{j}})_{+} ^{-\mathbf{\left( \frac{1}{2} + \frac{1-H}{q} \right)} } \right) \left( \prod _{j=1}^{r} (\mathbf{v}-{\mathbf{u}_{j})}_{+} ^{-\mathbf{\left( \frac{1}{2} + \frac{1-H}{q} \right)} }  \right)d\mathbf{v}\\
&=&c(\mathbf{H},q)  ^{2} \beta \left( \frac{1}{2}-\frac{1-{\bf H}}{q}, \frac{2-2{\bf H}}{q}\right)^{r}\int_{\mathbb{R} ^{d}} \int_{\mathbb{R} ^{d} } d\mathbf{u}d\mathbf{v}f({\bf u}) f({\bf v}) \vert {\bf u}- {\bf v}\vert ^{ \frac{2(H-{\bf 1})r}{q}} \\
&&\left( \prod _{j=1}^{q-r} (\mathbf{u}-\mathbf{y}_{j})_{+} ^{-\mathbf{\left( \frac{1}{2} + \frac{1-H}{q} \right)} }\right) \left( \prod _{j=q-r+1}^{2q-2r} (\mathbf{v}-\mathbf{y}_{j})_{+} ^{-\mathbf{\left( \frac{1}{2} + \frac{1-H}{q} \right)} } \right) 
\end{eqnarray*}
by using the Fubini theorem  and again relation  (13) in \cite{AT}, this leads to

\begin{eqnarray*}
&&\Vert Jf \otimes _{r} Jf \Vert ^{2} _{ L ^{2} ( \mathbb{R} ^{d (2q-2r)})}\\
&=&c(\mathbf{H},q)  ^{4} {\beta \left( \frac{1}{2}-\frac{1-{\bf H}}{q}, \frac{2-2{\bf H}}{q}\right)^{2r}\beta \left( \frac{1}{2}-\frac{1-{\bf H}}{q}, \frac{2-2{\bf H}}{q}\right)^{2q-2r}}\\
&&\int_{\mathbb{R} ^{d} } \int_{\mathbb{R} ^{d} } \int_{\mathbb{R} ^{d} } \int_{\mathbb{R} ^{d} }  d\mathbf{u}d\mathbf{v}d\mathbf{u}'d\mathbf{v}' ({\bf u}) f({\bf u} ') f({\bf v}) f({\bf v}') \\
&&\times \vert {\bf u}- {\bf v}\vert ^{ \frac{2(\bf{H}-{\bf 1})r}{q}}\vert {\bf u}'- {\bf v}'\vert ^{ \frac{2(\bf{H} -{\bf 1})r}{q}}\vert {\bf u}- {\bf u}'\vert ^{ \frac{2(\bf{H}-{\bf 1})(q-r)}{q}}
\vert {\bf v}- {\bf v}'\vert ^{ \frac{2(\bf{H}-{\bf 1})(q-r)}{q}}\\
&=&\frac{1}{q! ^{2} } ({\bf H} (2{\bf H}-{\bf 1})) ^{2}  \int_{\mathbb{R} ^{d} } \int_{\mathbb{R} ^{d} } \int_{\mathbb{R} ^{d} } \int_{\mathbb{R} ^{d} } d\mathbf{u}d\mathbf{v}d\mathbf{u}'d\mathbf{v}' f({\bf u}) f({\bf u} ') f({\bf v}) f({\bf v}') \\
&&\times \vert {\bf u}- {\bf v}\vert ^{ \frac{2(\bf H -{\bf 1})r}{q}}\vert {\bf u}'- {\bf v}'\vert ^{ \frac{2(\bf{H} -{\bf 1})r}{q}}\vert {\bf u}- {\bf u}'\vert ^{ \frac{2(\bf{H} -{\bf 1})(q-r)}{q}}
\vert {\bf v}- {\bf v}'\vert ^{ \frac{2(\bf{H}-{\bf 1})(q-r)}{q}}.
\end{eqnarray*}

The last quantity  converges to zero under assumption (\ref{c2-12}). \qed

Notice that $q=2$ and $d=1$ we retrieve the results in \cite{SlaTud2}. For $f=1$, the results in this section reduces to those in Theorem \ref{tat} from \cite{AT}.

\section{Applications to the stochastic heat equation with Hermite noise}

We will apply the main results in the previous section to some particular cases. First, we look to the solution to the heat equation driven by an Hermite noise. That is, we consider the following  linear stochastic heat equation driven by an additive  Hermite sheet  with $d+1$ parameters

\begin{equation}\label{heat}
\left\lbrace
\begin{array}{rcl}
 \frac{\partial u}{\partial t}(t,\mathbf{x})&=& \Delta u
(t,\mathbf{x}) +\dot Z_{H_{0}, {\bf H}} ^{q, d+1} (t,\mathbf{x}), \quad t\geq 0, \mathbf{x} \in \mathbb{R}^{d} \\ 
u(0,\mathbf{x})&=& 0, \quad \mathbf{x} \in \mathbb{R}^{d}
\end{array}
\right.
\end{equation}

We denoted by  $\Delta$  the Laplacian on $\mathbb{R} ^{d}$ and $Z^{q,d}_{H_{0}, {\bf H}}=\{Z^{q, d+1}_{H_{0}, \mathbf{H}}(t,\mathbf{x}); t \geq 0, \mathbf{x} \in \mathbb{R} ^{d}\}$ denotes  the $(d+1)$-parameter  Hermite sheet whose covariance is given by
\begin{equation*}
\mathbf{E} \left(  Z^{q, d+1}_{H_{0}, \mathbf{H}}(s,\mathbf{x})  Z^{q,d+1}_{H_{0}, \mathbf{H}}(t,\mathbf{y}) \right) =R_{H_{0}}(t,s) R_{\mathbf{H}} (\mathbf{x}, \mathbf{y} )
\end{equation*}
if $(H_{0},\mathbf{H} )= ({H_{0}}, H_{1},\ldots , H_{d})\in \left(\frac{1}{2}, 1\right) ^{d+1}$. We denoted by $\mathbf{H}= (H_{1},\ldots , H_{d})$ and
\begin{equation*}
R_{H}(t,s) = \frac{1}{2} ( \vert t\vert  ^{2H}+ \vert s\vert ^{2H}-\vert t-s\vert ^{2H}), \hskip0.3cm R_{\mathbf{H}} (\mathbf{x}, \mathbf{y} )= \prod_{j=1}^{d} R_{H_{j}} (x_{j}, y_{j})
\end{equation*}
if $s,t \in \mathbb{R}$ and ${\bf x}= (x_{1},.., x_{d}), {\bf y}=(y_{1},.., y_{d} ) \in \mathbb{R} ^{d}$.

The solution to (\ref{heat}) is understood in the mild sense. That is, the {\it mild} solution to (\ref{heat}) is a square-integrable process
$u=\{u(t,\mathbf{x}); t \geq 0, \mathbf{x} \in \mathbb{R} ^{d}\}$ defined by:
\begin{equation} \label{sol}
 u_{ H_{0}, {\bf H}}(t,\mathbf{x})=\int_{0}^{t}
\int_{\mathbb{R} ^{d}}G (t-s,\mathbf{x}-\mathbf{y})Z^{q,d+1}_{{H_{0}},\mathbf{H}}(ds,d\mathbf{y}), \hskip0.3cm t\geq 0, {\bf x} \in \mathbb{R} ^ {d}
\end{equation}
living in the space of jointly measurables random fields $\left( X(t,{\bf x}), t\geq 0, {\bf x} \in \mathbb{R}^{d} \right)$ such that for every $T>0$, $\sup_{t\in [0,T], {\bf x} \in \mathbb{R} ^{d}}\mathbf{E}\left| X(t, {\bf x})\right| ^{2} <\infty.$ 
 
The above integral is a Wiener integral with respect to the Hermite sheet, as introduced in Section 2 and $G(t,{\bf x}) $ is the Green function (or the fundamental solution) that satisfies $\frac{\partial u}{\partial t}-\Delta u=0$, i.e. 
\begin{equation}
\label{fund-sol-heat} 
G(t,\mathbf{x})=\left\{
\begin{array}{ll} (2 \pi t)^{-d/2} \exp\left( -\frac{|\mathbf{x}|^{2}}{2t}\right) & \mbox{if $t>0, \mathbf{x} \in \mathbb{R}^{d}$},  \\
0 & \mbox{if $t \leq 0, x \in \mathbb{R}^{d}$}.
\end{array} \right.
\end{equation}

The stochastic heat  equation (\ref{heat}) admits a unique mild solution $(u_{H_{0}, {\bf H }}(t,\mathbf{x}) )_{t\geq 0, \mathbf{x} \in \mathbb{R} ^{d} } $ if and only if  (see \cite{SlaTud1})

\begin{equation}
\label{cond1}d <4H_{0}+\sum_{i=1}^{d} (2H_{i}-1):=\gamma.
\end{equation}
In this case,  for every $T>0$,  $\underset{t\in [0,T], \mathbf{x}\in \mathbb{R} ^{d}} \sup   \mathbf{E}  \left(u(t, {\bf x})  ^{2}\right) < \infty. $\\

We will use   the following Parseval-type  formula (see Lemma A1 in \cite{BT}): for every $ f, g\in L ^ {2}(a,b)$ and  for every $0<\alpha <1$

\begin{equation}\label{parse}
\int_{a}^ {b}\int_{a} ^ {b}  du dv f(u)g(v) \vert u-v\vert ^ {-(1-\alpha)} = q_{\alpha} \int_{\mathbb{R}}\vert \tau \vert ^ {-\alpha} \mathcal{F}_{a,b} f(\tau) \overline{ \mathcal{F}_{a,b} g(\tau)}
\end{equation}
where  $(\mathcal{F} _{a,b} f)(\xi )= \int_{a}^{b} f(y) e ^{-i\xi y}dy$ (we use the notation $\mathcal{F}f=\mathcal{F} _{-\infty, \infty} f$) and
\begin{equation}\label{qa}
q_{\alpha} =(2 ^ {1-\alpha} \pi ^ {1/2}) ^ {-1} \frac{\Gamma (\alpha/2)}{\Gamma ((1-\alpha) /2)}.
\end{equation}
We recall that the Fourier transform of the function ${\bf y} \in \mathbb{R} ^{d} \to G(u, {\bf y})$ is $ \mathcal{F} G(u, \cdot) (\xi)= e^{-\frac{1}{2}u \vert \xi \vert ^{2}}.$

\subsection{Limit behavior of the solution when the Hurst index tends to  $1$}

The expression "Hurst index tends to  $1$"means that at least one component of the Hurst multi-index tends  to 1. We will apply Proposition \ref{pp1}  to obtain the asymptotic behavior of the solution (\ref{sol}) when  at least one of the  Hurst parameters $H_{0},H_{1},.., H_{d}$ converges to 1 and the other parameters are fixed.

\begin{theorem} 
Assume (\ref{cond1}) and let $A_{k}$ be as in (\ref{not}).  Fix $T>0$ and ${\bf x}\in \mathbb{R} ^{d}$.
Then
\begin{enumerate}

\item If

$$(H_{0}, {\bf H}_{A_{k}}) \to (1,.., 1) \in \mathbb{R} ^ {k+1} \mbox{ and }H_{j}, j\in \overline{A}_{k} \mbox{ are fixed} 
$$ 
then the stochastic process $(u_{H_{0}, {\bf H}}(t, {\bf x}), t\in [0,T])$ converges weakly  in $C[0,T]$  
to  the process $(u(t, {\bf x}), t\in [0,T])$ defined by 

\begin{equation}
\label{lim1}
u (t, {\bf x})=  \int_{0} ^ {t} du\int_{\mathbb{R} ^ {k} } d{\bf y}_{A_{k}} \int_{\mathbb{R} ^ {d-k} }d Z ^ {q, d-k} _{{\bf H}_{\overline{A}_{k}}} ({{\bf y}_{{\overline{A}_{k}}}} )G(t-u, {\bf x}-{\bf y}).
\end{equation}

\item If $ {\bf H}_{A_{k}} \to (1,.., 1) \in \mathbb{R} ^ {k}$ and $H_{0}, H_{j}, j\in \overline{A}_{k}$ are fixed, then $(u_{H_{0}, {\bf H}}(t, {\bf x}), t\in [0,T])$ converges weakly  in $C[0,T]$ to the stochastic process  $(u(t, {\bf x}), t\in [0,T])$ 
$$
u(t, {\bf x})= \int_{\mathbb{R} ^ {k} } d{\bf y}_{A_{k}}  \int_{0} ^ {t} \int_{\mathbb{R} ^ {d-k} }d Z ^ {q, d+1-k} _{H_{0},{\bf H}_{\overline{A}_{k}}} (u, {{\bf y}_{_{\overline{A}_{k}}}} )  G(t-u, {\bf x}-{\bf y}).
$$

\item If { $(H_{0}, {\bf H}) \to (1,...,1) \in \mathbb{R} ^ {d+1}$, then the weak  limit of $(u_{H_{0}, {\bf H}}(t, {\bf x}), t\in [0,T])$   in $C[0,T]$ is  $(u(t, {\bf x}), t\in [0,T])$  with

$$u(t, {\bf x})=
\left( \int_{0}^{t}
\int_{\mathbb{R} ^{d}}G (t-u,\mathbf{x}-\mathbf{y})d{\bf y}du \right) \frac{1}{\sqrt{q!}}H_{q}(Z). $$}
\end{enumerate}
\end{theorem}

\begin{remark}
As usual, by the weak convergence of the family  $(u_{H_{0}, {\bf H}}(t, {\bf x}), t\in [0,T])$ to $(u(t, {\bf x}), t\in [0,T])$ in $C[0, T]$ for fixed ${\bf x}\in \mathbb{R} ^{d}$ we mean the weak convergence of the  family of distributions of $u_{H_{0}, {\bf H}}(\cdot , {\bf x})$ to the law of $u(\cdot, {\bf x}) $ in $\left( C[0, T], \mathcal{B} (C[0,T])\right)$.

\end{remark}
\noindent {\bf Proof: }  Consider the function $F$ defined on $\mathbb{R}_{+}\times \mathbb{R}$ given by 
\begin{equation}
\label{27a-6}
F:  (u,{\bf y}) \to  1_{(0,t)} (u) (2\pi (t-u)) ^ {-\frac{d}{2}} e ^ {-\frac{\vert {\bf x}-{\bf y}\vert ^ {2}}{2(t-u)} }.
\end{equation}
We first show the convergence of finite dimensional distributions Consider the case 1. Let us show that this function belongs to $\vert \mathcal{H}_{H_{0}, {\bf H}}\vert\cap \mathcal{H}_{\overline{A}_{k}}$, with these two spaces defined by (\ref{normahb}) and (\ref{ha}) respectively. We know from \cite{BT} that, under (\ref{cond1}), the function $F$ (\ref{27a-6}) 
belongs to  the space $\vert \mathcal{H}_{H_{0}, {\bf H}}\vert .$

Let us  check that this function belongs to  the space $\mathcal{H}_{\overline{A}_{k}}$.  Writting

$$
F(u, {\bf y}) = F (u, {\bf y}_{A_{k}}, {\bf y}_{\overline{A}_{k}})= (2\pi u) ^{-\frac{d}{2}} e ^{-\frac{\vert {{\bf x}-\bf y}_{A_{k}}\vert ^{2}}{2u}}e ^{-\frac{\vert {{\bf x }- \bf y}_{\overline{A_{k}}}\vert ^{2}}{2u}} $$
we  have by the definition of the norm in  $\mathcal{H}_{\overline{A}}$ (see (\ref{ha})),

\begin{eqnarray*}
\Vert F\Vert _{\mathcal{H}_{\overline{A}_{k}} }&=&\sum_{j=1}^{k}  \int_{0}^{t} du \int_{\mathbb{R} ^{j}}d{\bf y} _{A_{j}} \left| \int_{\mathbb{R}^{d-j}} \int_{\mathbb{R} ^{d-j}} d{\bf y} _{\overline{A}_{j}}d{\bf z} _{\overline{A}_{j}}\right. \\
&&\left. \times  (2\pi u) ^{-\frac{d}{2}} e ^{-\frac{\vert {\bf y}_{A_{j}}\vert ^{2}}{2u}}e ^{-\frac{\vert {\bf y}_{\overline{A}_{j}}\vert ^{2}}{2u}} (2\pi u) ^{-\frac{d}{2}} e ^{-\frac{\vert {\bf y}_{A_{j}}\vert ^{2}}{2u}}e ^{-\frac{\vert {\bf z}_{\overline{A}_{j}}\vert ^{2}}{2u}} \vert {\bf y}_{\overline{A}_{j}}-{\bf z}_{\overline{A}_{j}}\vert ^{2H_{\overline{A}_{j}}-2} \right| ^{\frac{1}{2}}\\
&=&\sum_{j=1}^{k}\int_{0}^{t} du \int_{\mathbb{R} ^{j}}d{\bf y} _{A_{j}}  (2\pi u) ^{-\frac{j}{2}} e ^{-\frac{\vert {\bf y}_{A_{j}}\vert ^{2}}{2u}} \\
&&\times \left| \int_{\mathbb{R}^{d-j}} \int_{\mathbb{R} ^{d-j}} d{\bf y} _{\overline{A}_{j}}d{\bf z} _{\overline{A}_{j}}(2\pi u) ^{-\frac{d-j}{2}}e ^{-\frac{\vert {\bf y}_{\overline{A}_{j}}\vert ^{2}}{2u}} (2\pi u) ^{-\frac{d-j}{2}}e ^{-\frac{\vert {\bf z}_{\overline{A}_{j}}\vert ^{2}}{2u}} \vert {\bf y}_{\overline{A}_{j}}-{\bf z}_{\overline{A}_{j}}\vert ^{2H_{\overline{A}_{j}}-2} \right| ^{\frac{1}{2}}.
\end{eqnarray*}
By using Parseval's identity (\ref{parse})

$$ \int_{\mathbb{R}^{d-j}} \int_{\mathbb{R} ^{d-j}} d{\bf y} _{\overline{A}_{j}}d{\bf z} _{\overline{A}_{j}}(2\pi u) ^{-\frac{d-j}{2}}e ^{-\frac{\vert {\bf y}_{\overline{A}_{j}}\vert ^{2}}{2u}} (2\pi u) ^{-\frac{d-j}{2}}e ^{-\frac{\vert {\bf z}_{\overline{A}_{j}}\vert ^{2}}{2u}} \vert {\bf y}_{\overline{A}_{j}}-{\bf z}_{\overline{A}_{j}}\vert ^{2H_{\overline{A}_{j}}-2} =C_{j}\int_{\mathbb{R} ^{d-j}}d\xi e ^{-u\vert \xi \vert ^{2} } \vert \xi \vert ^{1-2{\bf H}_{\overline{A}_{j}}}$$
so with $C_{j}, C>0$
\begin{eqnarray*}
\Vert F\Vert _{\mathcal{H}_{\overline{A}_{k}} }&=& \sum_{j=1}^{k} C_{j}\int_{0}^{t} du \left| \int_{\mathbb{R} ^{d-j}}d\xi  e ^{-u\vert \xi \vert ^{2} } \vert \xi \vert ^{1-2{\bf H}_{\overline{A}_{j}}}\right| ^{\frac{1}{2}}\\
&=&C \int_{0}^{t} u ^{-\frac{d-j}{4} + \frac{1}{4} \sum _{a\in \overline{A}_{j}}(2H_{a}-1) }du
\end{eqnarray*}
and the last integral is finite if  for every $j=1,.., k$
\begin{equation}
\label{27a-5}
1-\frac{d-j}{4} + \frac{1}{4} \sum _{a\in \overline{A}_{j}}(2H_{a}-1)>0 \mbox{ or } d<4+j+ \sum _{a\in \overline{A}_{j}}(2H_{a}-1).
\end{equation}
The last bound is true due to (\ref{cond1}), so the function $F$ given by (\ref{27a-6}) belongs to $\vert \mathcal{H}_{H_{0}, {\bf H}}\vert\cap \mathcal{H}_{\overline{A}_{k}}$. 

Take $\lambda _{j}\in \mathbb{R}, t_{j}\geq 0$ for $j=1,..,N$ and denote by

\begin{equation}
\label{yn}
Y_{N}({\bf x})= \sum_{j=1}^ {N} \lambda _{j} u_{H_{0}, {\bf H}} (t_{j}, {\bf x}) = \int_{0} ^ {\infty} \int_{\mathbb{R} ^ {d} }  \left( \sum_{j=1}^ {N} \lambda _{j} 1_{(0, t_{j})} (u) G(t_{j}-u, {\bf x}-{\bf y})\right) dZ ^ {q, d+1} _{H_{0}, {\bf H}}(u, {\bf y}).
\end{equation}
From the above computations, the integrand $\sum_{j=1}^ {N} \lambda _{j} 1_{(0, t_{j})} (u) G(t_{j}-u, {\bf x}-{\bf y})$ in (\ref{yn}) belongs to $\vert \mathcal{H}_{H_{0}, {\bf H}}\vert\cap \mathcal{H}_{\overline{A}}$.  Therefore, by Proposition \ref{pp1},  the sequence $Y_{N}(x)$ (\ref{yn}) converges, as $(H_{0}, {\bf H}_{A_{k}}) \to (1,..,1) \in \mathbb{R} ^{k+1} $ to 
$$\sum_{j=1}^{N} \lambda _{j}\int_{0} ^ {t_{j}} du\int_{\mathbb{R} ^ {k} } d{\bf y}_{A_{k}} \int_{\mathbb{R} ^ {d-k} }d Z ^ {q, d-k} _{{\bf H}_{\overline{A}_{k}}} ({{\bf y}_{_{\overline{A}_{k}}}} )G(t_{j}-u, {\bf x}-{\bf y}) =\sum_{j=1}^{N} \lambda _{j}u(t_{j}, {\bf x})$$ 
with $u$ defined in (\ref{lim1}). This gives the convergence of the finite dimensional distribution of $(u_{H_{0}, {\bf H}}(t, {\bf x}), t\in [0,T])$ to the finite dimensional distributions of $(u(t, {\bf x}), t\in [0,T])$.

For the case 2., we have similarly
\begin{eqnarray*}
\Vert F\Vert _{\mathcal{H}_{\overline{A}_{k}} }&=& \sum_{j=1}^{k} C_{j} \left| \int_{0}^{t} \int_{0}^{t} dudv \vert u-v\vert ^{2H_{0}-2} \int_{\mathbb{R} ^{d-j}}d\xi  e ^{-\frac{1}{2}(u+v)\vert \xi \vert ^{2} } \vert \xi \vert ^{1-2{\bf H}_{\overline{A}_{j}}}\right| ^{\frac{1}{2}}\\
&=&C \left| \int_{0}^{t}\int_{0} ^{t} dudv \vert u-v\vert ^{2H_{0}-2}  (u+v) ^{-\frac{d-j}{2} + \frac{1}{2} \sum _{a\in \overline{A}_{j}}(2H_{a}-1) }\right| ^{\frac{1}{2}}
\end{eqnarray*}
and the above integral is finite under (\ref{cond1}). For the case 3., we notice in addition that the function $F$ given by (\ref{27a-6}) belongs to $ L ^{1} (\mathbb{R} ^{d+1})$.

Concerning the tightness, we recall that (see \cite{T}), for every $s,t\in [0,T], {\bf x} \in \mathbb{R} ^{d}$,

\begin{equation*}
\mathbf{E}\left| u_{H_{0}, {\bf H}} (t,{\bf x}) -u_{H_{0}, {\bf H}} (s,{\bf x})\right| ^{2} \leq C \vert t-s\vert ^{\gamma}
\end{equation*}
with $\gamma >0$ from (\ref{cond1}) and $C$ is a constant not depending on $s,t,{\bf x}$. Since $u_{H_{0}, {\bf H}} (t,{\bf x})$ is an element of the $(q+1)$th Wiener chaos, we use the hypercontractivity property for multiple stochastic integrals to get for every $p\geq 2$
\begin{equation}\label{tight}
\mathbf{E}\left| u_{H_{0}, {\bf H}} (t,{\bf x}) -u_{H_{0}, {\bf H}} (s,{\bf x})\right| ^{2p} \leq C \vert t-s\vert ^{\gamma p}
\end{equation}
and the tightness follows from (\ref{tight}) and the Billingsley criterion (see  \cite[Theorem 12.3]{B68} or  \cite{B99}). \qed 

\begin{remark}
Notice that when $(H_{0}, {\bf H}_{A_{k}}) \to (1,..,1) \in \mathbb{R} ^{k+1} $, the condition (\ref{cond1}) "converges" to  (\ref{27a-5}).

\end{remark}

\subsection{Limit behavior when the Hurst index tends to  $\frac{1}{2}$}

Fix $T>0$. When at least one of the components of Hurst multi-index goes to one-half, we have a central limit theorem. 

\begin{theorem}\label{tt3}
\begin{enumerate}
\item Assume 
\begin{equation}
\label{28a-3}
 (H_{0}, {\bf H}_{A_{k}}) \to \left( \frac{1}{2},...,\frac{1}{2}\right) \in \mathbb{R} ^{k+1}
\end{equation}
and
\begin{equation}
\label{28a-4}
d< 1+\frac{k}{2}+\sum_{a\in \overline{A}_{k}} H_{a}.
\end{equation}

Then the process ${(u_{H_{0}, {\bf H}}(t,{\bf x}), t\in [0, T]) }$ given by (\ref{sol}) converges weakly in $C[0,T]$ to the process $(u(t,{\bf x}), t\in [0,T] ) $ where $u$ is the mild solution to the heat equation 
\begin{equation}\label{heat1}
\left\lbrace
\begin{array}{rcl}
 \frac{\partial u}{\partial t}(t,\mathbf{x})&=& \Delta u
(t,\mathbf{x}) +\dot W_{H_{0}, {\bf H}} ^{q, d+1} (t,\mathbf{x}), \quad t>0, \mathbf{x} \in \mathbb{R}^{d} \\ 
u(0,\mathbf{x})&=& 0, \quad \mathbf{x} \in \mathbb{R}^{d}
\end{array}
\right.
\end{equation}
where $\left( W _{H_{0}, {\bf H}} (t, A_{1}\times A_{2} ), t\in [0, T], A_{1}\in \mathcal{B}_{b} (\mathbb{R} ^{k}), A_{2}\in  \mathcal{B}_{b} (\mathbb{R} ^{d-k})\right)$ is a Gaussian field with covariance
\begin{eqnarray*}
&&\mathbf{E}\left[ W_{H_{0}, {\bf H}}(t, {A_{1}}\times A_{2})W_{H_{0}, {\bf H}} (s, B_{1}\times B_{2})\right]\\
&=& (t\wedge s) \lambda _{k} (A_{1}\cap B_{1})  \int_{A_{2}\cap B_{2}} {\bf H}_{\overline{A}_{k}}(2{\bf H}_{\overline{A}_{k}}-{\bf 1})\vert {\bf y}_{\overline{A}_{k}}-{ {\bf z}_{\overline{A}_{k}}\vert}^{2{\bf H}_{\overline{A}_{k}}-2}d{\bf y}_{\overline{A}_{k}}d{\bf z}_{\overline{A}_{k}}.
\end{eqnarray*}
We denoted by $\lambda _{k}$ the Lebesque measure on $\mathbb{R}^{k}$. 

\item If $ {\bf H}_{A_{k}} \to \left( \frac{1}{2},...,\frac{1}{2}\right) \in \mathbb{R} ^{k}$ {, ${\bf H} _{B_{p}}\to
\left( 1,..,1 \right) \in \mathbb{R} ^{p}$}
and
\begin{equation}
\label{28a-9}
d< 2H+\frac{k}{2}+\sum_{a\in \overline{A}_{k}} H_{a}.
\end{equation}
then the process ${(u_{H_{0}, {\bf H}}(t,{\bf x}), t\in [0, T]) } $ given by (\ref{sol}) converges weakly in $C[0,T]$ to the process $(u(t,{\bf x}), t\in [0,T] ) $ where $u$ is the mild solution to the heat equation (\ref{heat1}) where the Gaussian noise has the following covariance 
\begin{eqnarray*}
&&\mathbf{E}\left[ W_{H_{0}, {\bf H}}(t, {A_{1}}\times A_{2})W_{H_{0}, {\bf H}} (s, B_{1}\times B_{2})\right]\\
&=& R_{H_{0}}(t,s) \lambda _{k} (A_{1}\cap B_{1})  \int_{A_{2}\cap B_{2}} {\bf H}_{\overline{A}_{k}}(2{\bf H}_{\overline{A}_{k}}-{\bf 1})\vert {\bf y}_{\overline{A}_{k}}- {{\bf z}_{\overline{A}_{k}}\vert}^{2{\bf H}_{\overline{A}_{k}}-2}d{\bf y}_{\overline{A}_{k}}d{\bf z}_{\overline{A}_{k}}.
\end{eqnarray*}

\item If $(H_{0}, {\bf H})\to \left( \frac{1}{2},...,\frac{1}{2}\right) \in \mathbb{R} ^{d+1}$ and $d=1$, then the weak limit of  $(u_{H_{0}, {\bf H}}, t\in [0, T]) $ in $C[0,T]$ is the solution to the heat equation (\ref{heat1}) driven by a space-time white noise. 
\end{enumerate}
\end{theorem}

\begin{remark}
The conditions (\ref{28a-4}), (\ref{28a-9}) and $d=1$ are the "limits" of (\ref{cond1}) in the cases 1., 2. and 3. respectively. 
\end{remark}
{\bf Proof: }  We will prove that the finite dimensional distributions of ${(u_{H_{0}, {\bf H}}(t,{\bf x}), t\in [0, T]) }  $ converge to those of $(u(t, {\bf x}), t\in [0,T])$ which satisfies (\ref{heat1}). In order to apply Proposition \ref{pp2}, we need to check conditions (\ref{c1-12}) and (\ref{c2-12}). 

\vskip0.2cm

\noindent {\it Checking condition (\ref{c1-12}). } Consider the case 1., i.e. assume (\ref{28a-3}) and (\ref{28a-4}).

Take $\lambda _{j}\in \mathbb{R}, t_{j}\geq 0$ for $j=1,..,N$ and denote by

$$Y_{N}({\bf x})= \sum_{j=1}^ {N} \lambda _{j} u_{H_{0}, {\bf H}} (t_{j}, {\bf x}) = \int_{0} ^ {\infty} \int_{\mathbb{R} ^ {d} }  \left( \sum_{j=1}^ {N} \lambda _{j} 1_{(0, t_{j})} (u) G(t_{j}-u, {\bf x}-{\bf y})\right) dZ ^ {q, d+1} _{H_{0}, {\bf H}}(u, {\bf y}).$$

We first check condition (\ref{c1-12}) for $ Y _{N}({\bf x})$. Let us calculate $\mathbf{E}\left( Y_{N} ({\bf x}  ) ^ {2}\right)$. By using the isometry (\ref{isoHW}), 

\begin{eqnarray*}
\mathbf{E}\left( Y_{N} ({\bf x}  ) \right)^ {2}&=&\sum_{j,k=1} ^ {N} \lambda _{j}\lambda _{k}H_{0}(2H_{0}-1){\bf H} (2{\bf H}-{\bf 1})\\
&& {\times} \int_{0}^ {t_{j}}du \int_{0}^ {t_{k}} dv \vert u-v \vert ^ {2H_{0}-2} \int_{\mathbb{R} ^ {d}}dy \int_{\mathbb{R}^ {d} }dz  G(t_{j}-u, {\bf x}-{\bf y})  G(t_{k}-v, {\bf x}-{\bf z}) \vert {\bf y}-{\bf z}\vert ^ {2{\bf H}-2}.
\end{eqnarray*}

\noindent Notice that, if ${\bf x}= (x ^ {(1)},.., x ^ {(d)}  ), {\bf y}= (y ^ {(1)},.., y ^ {(d)}  ),{\bf z}= (z ^ {(1)},.., z ^ {(d)}  )$ we have
$$G(t-u, {\bf x}-{\bf y}) =1_{(0,t)}(u) \prod_{a=1}^ {d} (2\pi (t-u) ) ^ {-\frac{d}{2}} e ^ {-\frac{\vert x^ {(a)}-y^ {(a)}\vert ^ {2}}{2(t-u) }}$$
and so
\begin{eqnarray*}
&& \int_{\mathbb{R} ^ {d}}dy \int_{\mathbb{R}^ {d} }dz  G(t_{j}-u, {\bf x}-{\bf y})  G(t_{k}-v, {\bf x}-{\bf z}) \vert {\bf y}-{\bf z}\vert ^ {2{\bf H}-2} \\
&=& \prod_{a=1} ^ {d}  \int_{\mathbb{R}} \int_{\mathbb{R}} dy^ {(a)} dz^ {(a)} (2\pi (t_{j}-u)) ^ {-\frac{1}{2}} (2\pi (t_{k}-v)) ^ {-\frac{1}{2}} e ^ {-\frac{\vert x^ {(a)}-y^ {(a)}\vert ^ {2}}{2(t_{j}-u) }} e ^ {-\frac{\vert x^ {(a)}-z^ {(a)}\vert ^ {2}}{2(t_{k}-v) }} \vert y ^ {(a)}-z^ {(a)} \vert ^ {2H_{a}-2}.
\end{eqnarray*}

We will apply the Parseval  identity  (\ref{parse})  with 
$$\alpha =2H_{a}-1  \mbox{ for every } a=1,..,d. $$
We get, for every $a=1,..,d$,
\begin{eqnarray*}
&& \int_{\mathbb{R}} \int_{\mathbb{R}} dy^ {(a)} dz^ {(a)} (2\pi (t_{j}-u)) ^ {-\frac{1}{2}} (2\pi (t_{k}-v)) ^ {-\frac{1}{2}} e ^ {-\frac{\vert x^ {(a)}-y^ {(a)}\vert ^ {2}} {2(t_{j}-u) }}  e ^ {-\frac{\vert x^ {(a)}-z^ {(a)}\vert ^ {2}}{2(t_{k}-v) }} \vert y ^ {(a)}-z^ {(a)} \vert ^ {2H_{a}-2}\\
&=& {q_{2H_{a}-1}} \int_{\mathbb{R}} d\tau \vert \tau \vert ^ {1-2H_{a}} e ^ {-\frac{1}{2} (t_{j}-u)\vert \tau \vert ^ {2}} e ^ {-\frac{1}{2} (t_{k}-v)\vert \tau \vert ^ {2}}.
\end{eqnarray*}

\noindent Now, by the change of variables $\tilde{\tau}=(t_{j}+t_{k}-2u) ^{\frac{1}{2}} \tau$,

\begin{eqnarray*}
 &&\int_{\mathbb{R}} d\tau \vert \tau \vert ^ {1-2H_{a}} e ^ {-\frac{1}{2} (t_{j}-u)\vert \tau \vert ^ {2}} e ^ {-\frac{1}{2} (t_{k}-v)\vert \tau \vert ^ {2}} \\
&=& (t_{j}+ t_{k}-u-v) ^ {-\frac{1}{2} +\frac{2H_{a}-1}{2}}\int_{\mathbb{R}} d\tau \vert \tau \vert ^ {1-2H_{a}} e ^ {-\frac{1}{2} \vert \tau \vert ^ {2}} =  (t_{j}+ t_{k}-u-v) ^ {H_{a}-1}\int_{\mathbb{R}} d\tau \vert \tau \vert ^ {1-2H_{a}} e ^ {-\frac{1}{2} \vert \tau \vert ^ {2}}. 
\end{eqnarray*}

\noindent Thus

\begin{eqnarray}
&&\mathbf{E}\left( Y_{N} ({\bf x}  ) \right)^ {2}\nonumber\\
&=&\sum_{j,k=1} ^ {N} \lambda _{j}\lambda _{k}H_{0}(2H_{0}-1){\bf H} (2{\bf H}-{\bf 1}) q_{2{\bf H}-{\bf 1}}\nonumber\\
&&\times\int_{0}^ {t_{j}}du \int_{0}^ {t_{k}} dv \vert u-v \vert ^ {2H_{0}-2}   (t_{j}+ t_{k} -u-v) ^{ H _{1}+...+ H_{d}-d} {\prod_{a=1}^{d}\int_{\mathbb{R}} d\tau \vert \tau \vert ^ {1-2H_{a}} e ^ {-\frac{1}{2} \vert \tau \vert ^ {2}}}\label{28a-6}
\end{eqnarray}
where ${q_{2H_{a}-1}}$ is defined in (\ref{qa}) and 
$$q_{2{\bf H}-{\bf 1}}= \prod_{a=1} ^{d} {q_{2H_{a}-1}}.$$

\noindent Notice that   for every $H\in (\frac{1}{2}, 1)$, we have 
$$ H(
2H-1) \Gamma(H-\frac{1}{2}) = H(2H-1) \frac{\Gamma (H+\frac{1}{2})}{H-\frac{1}{2}} \to _{H\to \frac{1}{2}} 2\Gamma (1) =2$$
and then
\begin{equation}
\label{28a-22}
H(2H-1) q_{2H-1} \to _{H\to \frac{1}{2}} (2\pi ) ^{-1}.
\end{equation}
Relation (\ref{28a-22}) implies
\begin{equation}
\label{28a-2}
{ \bf H} (2{\bf H}-{\bf 1}) q_{2{\bf H}-{\bf 1}} \to_{(H_{0}, {\bf H}_{A_{k}}) \to (\frac{1}{2},.., \frac{1}{2} ) \in \mathbb{R} ^{k+1} }(2\pi) ^{-k}q_{ 2{\bf H}_{\overline{A}_{k}}-{\bf 1}}.
\end{equation}
Let 
\begin{equation}
\label{ga}
\gamma:=H_{1}+...+H_{d}-d.
\end{equation}

\b We have, by integrating by parts 
\begin{eqnarray}
&&H_{0} (2H_{0}-1)  \int_{0} ^{t} \int_{0} ^{s} dudv \vert u-v\vert ^{2H_{0}-2}(t+s-u-v) ^{-\gamma}\nonumber 
\\
&=&H_{0} (2H_{0}-1)  \int_{0} ^{s} \int_{0} ^{s} dudv \vert u-v\vert ^{2H_{0}-2}(t+s-u-v) ^{-\gamma}\nonumber \\
&&+H_{0} (2H_{0}-1)  \int_{s} ^{t} \int_{0} ^{s} dudv \vert u-v\vert ^{2H_{0}-2}(t+s-u-v) ^{-\gamma}\nonumber\\
&=&H_{0} (2H_{0}-1)  2\int_{0} ^{s} \int_{0} ^{u} dudv \vert u-v\vert ^{2H_{0}-2}(t+s-u-v) ^{-\gamma}\nonumber\\
&&+H_{0} (2H_{0}-1)  \int_{s} ^{t} \int_{0} ^{s} dudv \vert u-v\vert ^{2H_{0}-2}(t+s-u-v) ^{-\gamma}\nonumber\\
&=&2H_{0}\int_{0} ^{s} du u ^{2H_{0}-1} (t+s-u) ^{-\gamma} \\
&&-H_{0} \int_{s}^{t}du u ^{2H_{0}-1}\left(  (t+s-u) ^{-\gamma} -(u-s) ^{2H_{0}-1} (t-u)^{-\gamma}\right)\nonumber \\
&&+2H_{0}\gamma \int_{0} ^{s}du \int_{0} ^{u}dv (u-v) ^{2H_{0}-1} (t+s-u-v) ^{-\gamma -1}\nonumber \\
&&+H_{0}\gamma  \int_{s} ^{t}du \int_{0} ^{u}dv (u-v) ^{2H_{0}-1} (t+s-u-v) ^{-\gamma -1}\nonumber\\
&=&2H_{0}\int_{0} ^{s} du u ^{2H_{0}-1} (t+s-u) ^{-\gamma} \nonumber\\
&&+H_{0} \int_{s}^{t}du u ^{2H_{0}-1}\left(  (t+s-u) ^{-\gamma } -(u-s) ^{2H_{0}-1} (t-u)^{-\gamma}\right)\nonumber \\
&&+ H_{0} \gamma \int_{0} ^{t}du \int_{0} ^{u}dv \vert u-v\vert  ^{2H_{0}-1} (t+s-u-v) ^{-\gamma -1}\label{28a-8}
\end{eqnarray}

\b Assuming (\ref{28a-3}), from (\ref{ga})
 $$\gamma\to -\left( d-\frac{k}{2}- \sum_{a\in \overline{A}_{k}} H_{a}\right):=\gamma_{0}$$
and, by taking the limit as $\gamma\to \gamma_{0}$ and $H_{0}\to \frac{1}{2}$ in (\ref{28a-8}), we get

\begin{eqnarray}
&&
H_{0} (2H_{0}-1)  \int_{0} ^{t} \int_{0} ^{s} dudv \vert u-v\vert ^{2H_{0}-2}(t+s-u-v) ^{-\gamma}\nonumber \\
&
\to&  \int_{0} ^{s} du  (t+s-u) ^{-\gamma_{0}}  \nonumber \\
&&+\frac{1}{2}\int_{s}^{t}du\left(  (t+s-u) ^{-\gamma_{0}} - (t-u)^{-\gamma_{0}}\right)\nonumber \\
&&+ \frac{1}{2}\gamma_{0}  \int_{0} ^{t}du \int_{0} ^{u}dv (t+s-u-v) ^{-\gamma_{0} -1}\nonumber \\
&=& \frac{1}{2} \frac{1} {(-\gamma_{0}+1)} \left( (t+s) ^{-\gamma_{0} +1} -\vert t-s\vert ^{-\gamma_{0} + 1}\right).\label{28a-7}
\end{eqnarray}
Consequently, as the limit (\ref{28a-3}) holds true, by plugging (\ref{28a-2}) and (\ref{28a-7}) into (\ref{28a-6}), we obtain

\begin{eqnarray*}
\mathbf{E} Y_{N}({\bf x}) ^{2} &\to& \frac{1}{2} \frac{1} {-\gamma_{0}+1} (2\pi) ^{-k} \sum_{j,k=1}^{N} \lambda _{j} \lambda _{k} \left(  (t_{j}+ t_{k} ) ^{-\gamma_{0} +1}- \vert t_{j}-t_{k} \vert ^{-\gamma_{0}+ 1}\right) q_{2{\bf H}_{\overline{ A}_{k}}-{\bf 1}}\\
&&\times  \prod_{a\in A_{k}} \int_{\mathbb{R}} d\tau e ^ {-\frac{1}{2} \vert \tau \vert ^ {2}}
\prod_{a\in \overline{A}_{k} } \int_{\mathbb{R}} d\tau \vert \tau \vert ^ {1-2H_{a}} e ^ {-\frac{1}{2} \vert \tau \vert ^ {2}}\\
&=& \frac{1}{2} \frac{1} {-\gamma_{0}+1}  (2\pi) ^{-k} \sum_{j,k=1}^{N} \lambda _{j} \lambda _{k} \left(  (t_{j}+ t_{k} ) ^{-\gamma_{0} +1}- \vert t_{j}-t_{k} \vert ^{-\gamma_{0}+ 1}\right) \\
&&\times q_{2{\bf H}_{\overline{ A}_{k}}-{\bf 1}} (\sqrt{2\pi} ) ^{k} \prod_{a\in \overline{A}_{k} } \int_{\mathbb{R}} d\tau \vert \tau \vert ^ {1-2H_{a}} e ^ {-\frac{1}{2} \vert \tau \vert ^ {2}}\\
&=&\frac{1}{2} \frac{1} {-\gamma_{0}+1}  (2\pi) ^{-\frac{k}{2}} \sum_{j,k=1}^{N} \lambda _{j} \lambda _{k} \left(  t_{j}+ t_{k} ) ^{-\gamma_{0} +1}- \vert t_{j}-t_{k} \vert ^{-\gamma_{0}+ 1}\right) \\
&&\times q_{2{\bf H}_{\overline{ A}_{k}}-{\bf 1}} \prod_{a\in \overline{A}_{k} } \int_{\mathbb{R}} d\tau \vert \tau \vert ^ {1-2H_{a}} e ^ {-\frac{1}{2} \vert \tau \vert ^ {2}}.
\end{eqnarray*}

\b On the other hand, if $u$ is the solution to (\ref{heat1}), then

\begin{eqnarray*}
 &&\mathbf{E}\left(\sum_{j=1} ^{N} \lambda _{j}  u(t_{j}, {\bf x} )\right)^{2} = \sum_{j,k=1} ^{N} \lambda _{j} \lambda _{k} \int_{0} ^{t_{j}\wedge t_{k}}  du \int_{\mathbb{R} ^{k} } d{\bf y} _{A_{k}} \int_{\mathbb{R} ^{d-k}} \int_{\mathbb{R} ^{d-k}} d{\bf y}_{\overline{A}_{k}}d{\bf z}_{\overline{A}_{k}}\\
&&\times  (2\pi (t_{j}-u) ^{-\frac{d}{2}} e ^{-\frac{ \vert {\bf y} _{A_{k}}\vert ^{2} }{2(t_{j}-u)}}e ^{-\frac{ \vert {\bf y} _{\overline{A}_{k}}\vert ^{2} }{2(t_{j}-u)}} (2\pi (t_{k}-u) ^{-\frac{d}{2}} e ^{-\frac{ \vert {\bf y} _{A_{k}}\vert ^{2} }{2(t_{k}-u)}}e ^{-\frac{ \vert {\bf z} _{\overline{A}_{k}}\vert ^{2} }{2(t_{k}-u)}} \\
&=&\sum_{j,k=1} ^{N} \lambda _{j} \lambda _{k}\int_{0} ^{t_{j}\wedge t_{k}}  du (2\pi) ^{-k}  \int_{\mathbb{R}^{k} } d\xi e ^{-(t_{j}+ t_{k}-2u)\vert \xi \vert ^{2}} q_{2{\bf H}_{\overline{A}_{k}}-{\bf 1}}\int_{\mathbb{R}^{d-k}} d\tau e ^{-(t_{j}+ t_{k}-2u)\vert \tau \vert ^{2}} \vert \tau \vert ^{\frac{1}{2}\sum_{a\in \overline{A}_{k}} (2H_{a}-1)}\\
&=&\sum_{j,k=1} ^{N} \lambda _{j} \lambda _{k} \int_{0} ^{t_{j}\wedge t_{k}} du    (t_{j}+t_{k}-2u) ^{-\gamma _{0} } (2\pi) ^{-k} \int_{\mathbb{R}^{k} } d\xi e ^{-\vert \xi \vert ^{2}} q_{2{\bf H}_{\overline{A}_{k}}-{\bf 1}} \prod_{a\in {\overline{A}_{k}}} \int_{\mathbb{R}} d\tau e ^{-\vert \tau \vert ^{2}} \vert \tau \vert ^{1-2H_{a}}\\
&=& \frac{1}{2} \frac{1}{-\gamma_{0}+1} (2\pi) ^{-\frac{k}{2}}\sum_{j,k=1} ^{N} \lambda _{j} \lambda _{k}  
\left( (t_{j}+t_{k}) ^{-\gamma_{0}+1} -\vert t_{j}-t_{k} \vert ^{-\gamma_{0}+1} \right) q_{2{\bf H}_{\overline{A}_{k}}-{\bf 1}} \prod_{a\in {\overline{A}_{k}}} \int_{\mathbb{R}} d\tau e ^{-\vert \tau \vert ^{2}} \vert \tau \vert ^{1-2H_{a}}.
\end{eqnarray*}

\b The point 2. follows similarly. Let us discuss point 3.  Assume $H_{0},H_{1},.., H_{d} $ converge all to $\frac{1}{2}$. Notice that in this case condition (\ref{cond1}) implies $d<2$ so $d=1$! Then, from (\ref{ga})

$$\gamma \to \frac{d}{2}=\frac{1}{2}.$$

\b Therefore, from (\ref{28a-8}), as $H_{0},H_{1},.., H_{d}\to \frac{1}{2}$
\begin{eqnarray}
&& H_{0} (2H_{0}-1)  \int_{0} ^{t} \int_{0} ^{s} dudv \vert u-v\vert ^{2H_{0}-2}(t+s-u-v) ^{-\gamma}\nonumber \\
&&\to  \frac{1}{2} \int_{0} ^{t}du \left( (t-u) ^{-\frac{1}{2}} -(t+s-u) ^{-\frac{1}{2}}\right)
+ \frac{1}{2} \times\frac{1}{2} \int_{0}^{t}du \int_{0} ^{s}dv (t+s-u-v) ^{-\frac{3}{2}}\nonumber \\
&=& \left( (t+s) ^{\frac{1}{2} }- \vert t-s\vert ^{\frac{1}{2}}\right).\label{28a-10}
\end {eqnarray}
and we obtain, by combining (\ref{28a-10}) and (\ref{28a-6}), by taking the limit (\ref{28a-3})
\begin{eqnarray*}
 \mathbf{E} Y_{N}({\bf x}) ^{2} &\to& (2\pi) ^{-1} \sum_{j,k=1}^{N} \lambda _{j} \lambda _{k} \left(  (t_{j}+ t_{k} ) ^{\frac{1}{2}} -\vert t_{j}-t_{k} \vert ^{\frac{1}{2}}\right) \int_{\mathbb{R}} d\tau e ^{-\frac{1}{2} \vert \tau \vert ^{2}}\\
&=&\sum_{j,k=1}^{N} \lambda _{j} \lambda _{k} \left(  (t_{j}+ t_{k} ) ^{\frac{1}{2}} -\vert t_{j}-t_{k} \vert ^{\frac{1}{2}}\right)  \sqrt{2\pi} \\
&=&(2\pi) ^{-\frac{1}{2}} \sum_{j,k=1}^{N} \lambda _{j} \lambda _{k} \left(  (t_{j}+ t_{k} ) ^{\frac{1}{2}} -\vert t_{j}-t_{k} \vert ^{\frac{1}{2}}\right)\end{eqnarray*}
which coincides with the $ \mathbf{E}\left( \sum_{j=1}^{N} \lambda _{j} u(t_{j}, {\bf x})\right) ^{2} $ where $u$ is the solution of the heat equation (\ref{heat1}) driven  by a space-time white noise (see \cite{Sw} or \cite{T}).

\vskip0.2cm

{\it Checking condition  (\ref{c2-12}). } In order to check condition (\ref{c2-12}), we need to show in the case 1. (the other situations are similar) that for every $t_{1}, t_{2}, t_{3}, t_{4} \in [0, T]$,
\begin{eqnarray*}
I:&=&\sup_{(H_{0}, {\bf H }_{A_{k}})\in [\frac{1}{2}, 1] ^{k+1}}   \int_{0} ^{t_{1}} du_{1}...\int_{0} ^{t_{4}} du_{4} \vert u_{1}-u_{2} \vert ^{-\alpha_{0}}   \vert u_{2}-u_{3} \vert ^{-\alpha_{0}}   \vert u_{3}-u_{4} \vert ^{-\beta_{0}}   \vert u_{4}-u_{1} \vert ^{-\beta_{0}}\\
 &&\times  \int_{\mathbb{R } ^{d}} d{\bf y}_{1}...\int_{\mathbb{R } ^{d}} d{\bf y}_{4} \frac{1}{(2\pi (t_{1}-u_{1}))^{{\frac{d}{2}}}} e ^{-\frac{ {\vert{\bf x}-{\bf y}_{1}\vert^{2}} }{2(t_{1}-u_{1})}} \frac{1}{(2\pi (t_{2}-u_{2}))^{{\frac{d}{2}}}} e ^{-\frac{ {\vert{\bf x}-{\bf y}_{2}\vert^{2}} }{2(t_{2}-u_{2})}} \\
&&\times\frac{1}{(2\pi (t_{3}-u_{3}))^{{\frac{d}{2}}}} e ^{-\frac{ {\vert{\bf x}-{\bf y}_{3}\vert^{2} }}{2(t_{3}-u_{3})}} \frac{1}{(2\pi (t_{4}-u_{4}))^{{\frac{d}{2}}}} e ^{-\frac{ {\vert{\bf x}-{\bf y}_{4}\vert^{2}}}{2(t_{4}-u_{4})}}\\
&&\vert {\bf y}_{1}-{\bf y}_{2}\vert ^{-\alpha} \vert {\bf y}_{2}-{\bf y}_{3}\vert ^{-\alpha} \vert {\bf y}_{3}-{\bf y}_{4}\vert ^{-\beta} \vert {\bf y}_{4}-{\bf y}_{1}\vert ^{-\beta }<\infty 
\end{eqnarray*}
 with
$$\alpha = \frac{2(1-{\bf H}) r}{q}, \hskip0.2cm \beta =\frac{2(1-{\bf H }) (q-r) }{q} , \hskip0.2cm \alpha _{0} =\frac{2(1-H_{0})r}{q}, \hskip0.2cm \beta _{0}= \frac{2(1-H_{0} )(q-r)}{q}$$
for every $r=1,.., q-1.$ After the change of variables $t_{i}-u_{i}= \tilde{u}_{i}, \tilde{{\bf y} } = {\bf x}- {\bf y}$, we will have to show that
\begin{eqnarray*}
I&=&\sup_{(H_{0}, {\bf H }_{A_{k}})\in [\frac{1}{2}, 1] ^{k+1}}   \int_{0} ^{t_{1}} du_{1}...\int_{0} ^{t_{4}} du_{4}\\
&& \vert u_{1}-u_{2} -(t_{1}-t_{2} )\vert ^{-\alpha_{0}}   \vert u_{2}-u_{3} -(t_{2}-t_{3}) \vert ^{-\alpha_{0}}   \vert u_{3}-u_{4} -(t_{3}-t_{4})\vert ^{-\beta_{0}}   \vert u_{4}-u_{1} -(t_{4}-t_{1})\vert ^{-\beta_{0}}  \\
&&\int_{\mathbb{R } ^{d}} d{\bf y}_{1}...\int_{\mathbb{R } ^{d}} d{\bf y}_{4} \frac{1}{(2\pi u_{1})^{{\frac{d}{2}}}} e ^{-\frac{ {-\vert{\bf y}_{1}\vert}^{2}}{2u_{1}}} \frac{1}{(2\pi u_{2})^{{\frac{d}{2}}}} e ^{-\frac{ {-\vert{\bf y}_{2}\vert}^{2}}{2u_{2}}} \frac{1}{(2\pi u_{3})^{{\frac{d}{2}}}} e ^{-\frac{ {-\vert{\bf y}_{3}\vert}^{2}}{2u_{3}}} \frac{1}{(2\pi u_{4})^{{\frac{d}{2}}}}e ^{-\frac{ {-\vert{\bf y}_{4}\vert}^{2}}{2u_{4}}} \\
&&\vert {\bf y}_{1}-{\bf y}_{2}\vert ^{-\alpha} \vert {\bf y}_{2}-{\bf y}_{3}\vert ^{-\alpha} \vert {\bf y}_{3}-{\bf y}_{4}\vert ^{-\beta} \vert {\bf y}_{4}-{\bf y}_{1}\vert ^{-\beta }<\infty. 
\end{eqnarray*}
Next, we write for the integrals $d{\bf y}_{i}$
\begin{eqnarray*}
&&\int_{\mathbb{R } ^{d}} d{\bf y}_{1}...\int_{\mathbb{R } ^{d}} d{\bf y}_{4} ...\\
&=& \prod_{j\in A_{k}} \int_{\mathbb{R}} d y_{1} ^{(j)} ...\int_{\mathbb{R}} d y_{4} ^{(j)} \frac{1}{\sqrt{2\pi u_{1}}} e ^{-\frac{ {\vert y_{1} ^{(j)} \vert ^{2}}}{2u_{1}}} ...\frac{1}{\sqrt{2\pi u_{4}}} e ^{-\frac{ {\vert y_{4} ^{(j)} \vert ^{2}}}{2u_{4}}} \\
&&\times \prod_{j\in \overline{A}_{k}} \int_{\mathbb{R}} d y_{1} ^{(j)} ...\int_{\mathbb{R}} d y_{4} ^{(j)} \frac{1}{\sqrt{2\pi u_{1}}} e ^{-\frac{{\vert y_{1} ^{(j)} \vert  ^{2}}}{2u_{1}}} ...\frac{1}{\sqrt{2\pi u_{4}}} e ^{-\frac{{\vert y_{4} ^{(j)} \vert ^{2}}}{2u_{4}}}.
\end{eqnarray*}
We will separate the integral $dy_{1}^{(j)} $, for every $j=1,..,d$,  as follows 
$$ \int_{\mathbb{R}} d y_{1} ^{(j)}= \int_{ \vert y_{1}\vert ^{(j)} >\sqrt{2T}} dy_{1} ^{(j)} + \int_{ {\vert y_{1}\vert ^{(j)} \leqslant\sqrt{2T}}} dy_{1} ^{(j)}$$
and similarly for the integrals $dy_{2} ^{(j)}, dy_{3} ^{(j)}, dy_{4} ^{(j)}.$  We use the fact that on the set 
$$y^{2} >2T>2u$$
the function 
$$u\to \frac{1}{\sqrt{u} } e ^{-\frac{y ^{2}}{2u}}  \mbox{ is increasing }$$
and we majorize
$$\frac{1}{\sqrt{u} } e ^{-\frac{y ^{2}}{2u}}  \mbox{ by } \frac{1}{\sqrt{T} } e ^{-\frac{y ^{2}}{2T}} $$
On the other hand, on the set 
$$y^{2} {\leqslant }2T$$
we majorize 
$$\frac{1}{\sqrt{u} } e ^{-\frac{y ^{2}}{2u}}  \mbox{ by  a constant. }  $$ 

\b In this way, the quantity $I$ can be bounded by 
\begin{eqnarray*}
I&\leq &C\sup_{(H_{0}, {\bf H }_{A_{k}})\in [\frac{1}{2}, 1] ^{k+1}}   \int_{0} ^{t_{1}} du_{1}...\int_{0} ^{t_{4}} du_{4}\\
&& \vert u_{1}-u_{2} -(t_{1}-t_{2} )\vert ^{-\alpha_{0}}   \vert u_{2}-u_{3} -(t_{2}-t_{3}) \vert ^{-\alpha_{0}}   \vert u_{3}-u_{4} -(t_{3}-t_{4})\vert ^{-\beta_{0}}   \vert u_{4}-u_{1} -(t_{4}-t_{1})\vert ^{-\beta_{0}}  \\
&&\prod_{j\in A_{k}} \int_{\mathbb{R}} d y_{1} ^{(j)} ...\int_{\mathbb{R}} d y_{4} ^{(j)}\left(  \frac{1}{\sqrt{2\pi T}} e ^{-{\frac{ \vert y^{(j)}_{1}\vert ^{2}}{2T}}}1_{ \vert y_{1} ^{(j)}\vert  >\sqrt{2T}} + 1_{ {\vert y_{1}\vert ^{(j)} \leqslant\sqrt{2T}}} \right)\\
&&\times...\left(  \frac{1}{\sqrt{2\pi T}} e ^{-{\frac{ \vert y^{(j)}_{4}\vert ^{2}}{2T}}}1_{ \vert y_{4} ^{(j)}\vert  >\sqrt{2T}} + 1_{{\vert y_{4}\vert ^{(j)} \leqslant\sqrt{2T}}} \right)\\
&&\times \vert y_{1} ^{(j)} -y_{2} ^{(j)} \vert ^{-\alpha_{j}} \vert y_{2} ^{(j)} -y_{3} ^{(j)} \vert ^{-\alpha_{j}}\vert y_{3} ^{(j)} -y_{4} ^{(j)} \vert ^{-\beta_{j}}\vert y_{4} ^{(j)} -y_{1} ^{(j)} \vert ^{-\beta_{j}}\\
&&\times R 
\end{eqnarray*}
with $\alpha_{j}= \frac{2(1-H_{j})r}{q}, \beta_{j}=\frac{2(1-H_{j})(q-r)}{q}$ for every $j=1,..,d$ and 
\begin{eqnarray*}
R&=&\prod _{j\in \overline{A}_{k}}\int_{\mathbb{R}} d y_{1} ^{(j)} ...\int_{\mathbb{R}} d y_{4} ^{(j)}\left(  \frac{1}{\sqrt{2\pi T}} e ^{-{\frac{ \vert y^{(j)}_{1}\vert ^{2}}{2T}}}1_{ \vert y_{1} ^{(j)} \vert \geqslant\sqrt{2T}} + 1_{ {\vert y_{1}\vert ^{(j)} \leqslant\sqrt{2T}}} \right)\\
&&\times...\left(  \frac{1}{\sqrt{2\pi T}} e ^{-{\frac{ \vert y^{(j)}_{4}\vert ^{2}}{2T}}}1_{ \vert y_{4} ^{(j)} \vert \geqslant\sqrt{2T}} + 1_{ {\vert y_{4}\vert ^{(j)} \leqslant\sqrt{2T}}} \right)\\
&&\times \vert y_{1} ^{(j)} -y_{2} ^{(j)} \vert ^{-\alpha_{j}} \vert y_{2} ^{(j)} -y_{3} ^{(j)} \vert ^{-\alpha_{j}}\vert y_{3} ^{(j)} -y_{4} ^{(j)} \vert ^{-\beta_{j}}\vert y_{4} ^{(j)} -y_{1} ^{(j)} \vert ^{-\beta_{j}}.
\end{eqnarray*}
Consequently, we can write
\begin{eqnarray*}
I&\leq & C \sup_{H_{0} \in [\frac{1}{2}, 1]}  \int_{0} ^{t_{1}} du_{1}...\int_{0} ^{t_{4}} du_{4}\\
&& \vert u_{1}-u_{2} -(t_{1}-t_{2} )\vert ^{-\alpha_{0}}   \vert u_{2}-u_{3} -(t_{2}-t_{3}) \vert ^{-\alpha_{0}}   \vert u_{3}-u_{4} -(t_{3}-t_{4})\vert ^{-\beta_{0}}   \vert u_{4}-u_{1} -(t_{4}-t_{1})\vert ^{-\beta_{0}}  \\
&&\sup _{ {\bf H }_{A_{k}}\in [\frac{1}{2}, 1] ^{k}}\prod_{j\in A_{k}} \int_{\mathbb{R}} d y_{1} ^{(j)} ...\int_{\mathbb{R}} d y_{4} ^{(j)}\left(  \frac{1}{\sqrt{2\pi T}} e ^{-{\frac{ \vert y^{(j)}_{1}\vert ^{2}}{2T}}}1_{ \vert y_{1}^{(j)} \vert \geqslant \sqrt{2T}} + 1_{{\vert y^{(j)}_{1} \vert \leqslant\sqrt{2T}}} \right)\\
&& .....\left(  \frac{1}{\sqrt{2\pi T}} e ^{-{\frac{ \vert y^{(j)}_{4}\vert ^{2}}{2T}}}1_{ \vert y_{4} ^{(j)} \vert \geqslant\sqrt{2T}} + 1_{ {\vert y_{4}^{(j)}\vert  \leqslant \sqrt{2T}}} \right)\\
&&\vert y_{1} ^{(j)} -y_{2} ^{(j)} \vert ^{-\alpha_{j}} \vert y_{2} ^{(j)} -y_{3} ^{(j)} \vert ^{-\alpha_{j}}\vert y_{3} ^{(j)} -y_{4} ^{(j)} \vert ^{-\beta_{j}}\vert y_{4} ^{(j)} -y_{1} ^{(j)} \vert ^{-\beta_{j}}\\
&&\times R.
\end{eqnarray*}
Note that $R$ does not depend on $H_{0}, {\bf H}_{A_{k}}$ and
 \begin {eqnarray*}
 &&\sup_{H_{0} \in [\frac{1}{2}, 1]}  \int_{0} ^{t_{1}} du_{1}...\int_{0} ^{t_{4}} du_{4}\\
&& \vert u_{1}-u_{2} -(t_{1}-t_{2} )\vert ^{-\alpha_{0}}   \vert u_{2}-u_{3} -(t_{2}-t_{3}) \vert ^{-\alpha_{0}}   \vert u_{3}-u_{4} -(t_{3}-t_{4})\vert ^{-\beta_{0}}   \vert u_{4}-u_{1} -(t_{4}-t_{1})\vert ^{-\beta_{0}}  \\
&\leq & \int_{0} ^{T}du_{1}...\int_{0} ^{T} du_{4}\vert u_{1}-u_{2} \vert ^{-\alpha_{0}}   \vert u_{2}-u_{3} \vert ^{-\alpha_{0}}   \vert u_{3}-u_{4} \vert ^{-\beta_{0}}   \vert u_{4}-u_{1} \vert ^{-\beta_{0}} 
\end{eqnarray*}
which is finite by Lemma 3.3 in \cite{BaiTa} since
$$2\alpha +2\beta + 4 = 2(2H-2)+4 = 4H>1.$$
Therefore, in order to conclude, it remains to show that 
\begin{eqnarray*}
&&\sup _{ {\bf H }_{A_{k}}\in [\frac{1}{2}, 1] ^{k}}\prod_{j\in A_{k}} \int_{\mathbb{R}} d y_{1} ^{(j)} ...\int_{\mathbb{R}} d y_{4} ^{(j)}\\
&&\left(  \frac{1}{\sqrt{2\pi T}} e ^{-{\frac{ \vert y^{(j)}_{1}\vert ^{2}}{2T}}}1_{ \vert y_{1} ^{(j)} \vert \geqslant \sqrt{2T}} + 1_{ {\vert y_{1}^{(j)}\vert } \leqslant\sqrt{2T}} \right)...\left(  \frac{1}{\sqrt{2\pi T}} e ^{-{\frac{ \vert y^{(j)}_{4}\vert ^{2}}{2T}}}1_{ \vert y_{4} ^{(j)} \vert \geqslant \sqrt{2T}} + 1_{ {\vert y_{4}^{(j)}\vert } \leqslant \sqrt{2T}} \right) \\
&&\times \vert y_{1} ^{(j)} -y_{2} ^{(j)} \vert ^{-\alpha_{j}} \vert y_{2} ^{(j)} -y_{3} ^{(j)} \vert ^{-\alpha_{j}}\vert y_{3} ^{(j)} -y_{4} ^{(j)} \vert ^{-\beta_{j}}\vert y_{4} ^{(j)} -y_{1} ^{(j)} \vert ^{-\beta_{j}}<\infty.
\end{eqnarray*}

Assume for simplicity $A_{k}=\{1,2,.., k\}. $  To check that the above quantity is finite, it suffices to prove that 

\begin{eqnarray*}
&&\sup_{H\in [\frac{1}{2}, 1]} \int_{\mathbb{R} } dy_{1}....\int_{\mathbb{R} } dy_{4} \left( e ^{-{\frac{ \vert y_{1}\vert ^{2}}{2T}}}1_{ \vert y_{1} \vert \geqslant \sqrt{2T}} + 1_{ {\vert y_{1}\vert } \leqslant \sqrt{2T}} \right)....\left( e ^{-{\frac{ \vert y_{4}\vert ^{2}}{2T}}}1_{ \vert y_{4}\vert  \geqslant \sqrt{2T}} + 1_{{\vert y_{4}\vert } \leqslant\sqrt{2T}} \right)\\
&&\times \vert y_{1}-y_{2} \vert ^{-\alpha } \vert y_{2}-y_{3} \vert ^{-\alpha } \vert y_{3}-y_{4} \vert ^{-\beta} \vert y_{4}-y_{1} \vert ^{-\beta } <\infty. 
\end{eqnarray*}

Using $\prod_{i=1}^{4} (A_{i}+B_{i})= A_{1}A_{2}A_{3}A_{4}+A_{1}B_{2}B_{3}B_{4}+....+B_{1}B_{2} B_{3}B_{4}$, the last integrals can be expressed as a sum of several terms, involving integrals on the sets {$\vert y_{i}\vert  \geqslant \sqrt{2T}$ and  $\vert y_{i} \vert \leqslant \sqrt{2T}$}.

Let us start with the first summand, namely

\begin{eqnarray*}
T_{1}&:=&\sup_{H\in [\frac{1}{2}, 1]} \int_{\mathbb{R} } dy_{1}....\int_{\mathbb{R} } dy_{4} e ^{-\frac{ y_{1}   ^{2}}{2T}}1_{ \vert y_{1}\vert \geqslant  \sqrt{2T}} e ^{-\frac{ y_{2}   ^{2}}{2T}}1_{ \vert y_{2} \vert \geqslant  \sqrt{2T}} e ^{-\frac{ y_{3}   ^{2}}{2T}}1_{ \vert y_{3} \vert \geqslant  \sqrt{2T}} e ^{-\frac{ y_{4}   ^{2}}{2T}}1_{ \vert y_{4} \vert \geqslant \sqrt{2T}}\\
&&\times \vert y_{1}-y_{2} \vert ^{-\alpha } \vert y_{2}-y_{3} \vert ^{-\alpha } \vert y_{3}-y_{4} \vert ^{-\beta} \vert y_{4}-y_{1} \vert ^{-\beta } . 
\end{eqnarray*}
Since $\vert y_{1}-y_{2}\vert ^{2} \leq 2(y_{1} ^{2} + y_{2} ^{2} ) $ we have

\begin{equation}
\label{28a-11}
 \vert y_{1}-y_{2}\vert ^{2}+\vert y_{2}-y_{3}\vert ^{2}+\vert y_{3}-y_{4}\vert ^{2}+\vert y_{4}-y_{1}\vert ^{2}\leq 4 (y_{1} ^{2} + y_{2} ^{2}+y_{3} ^{2} + y_{4} ^{2})
\end{equation}
so 

$$e ^{-\frac{ y_{1}   ^{2} + y_{2} ^{2} + y_{3} ^{2} + y_{4} ^{2} }{2T}}\leq e ^{-\frac{1}{8T} (\vert y_{1}-y_{2}\vert ^{2} +\vert y_{2}-y_{3}\vert ^{2} +\vert y_{3}-y_{4}\vert ^{2} +\vert y_{4}-y_{1}\vert ^{2})}.$$

\b Hence,  $T_{1}$ can be  bounded as follows

\begin{eqnarray*}
T_{1}&\leq &  \sup_{H\in [\frac{1}{2}, 1]} \int_{\mathbb{R} } dy_{1}....\int_{\mathbb{R} } dy_{4} e ^{-\frac{1}{8T} (\vert y_{1}-y_{2}\vert ^{2} +\vert y_{2}-y_{3}\vert ^{2} +\vert y_{3}-y_{4}\vert ^{2} +\vert y_{4}-y_{1}\vert ^{2})}\\
&&\vert y_{1}-y_{2} \vert ^{-\alpha } \vert y_{2}-y_{3} \vert ^{-\alpha } \vert y_{3}-y_{4} \vert ^{-\beta} \vert y_{4}-y_{1} \vert ^{-\beta } . 
\end{eqnarray*}

\b We apply the power counting theorem, see the Appendix. Consider the set of affine functionals 
$$T'=\{y_{1}-y_{2}, y_{2}-y_{3}, y_{3}-y_{4}, y_{4}-y_{1}\}. $$ 
The only padded subset of $T'$ is $T'$ itself.  We apply the  power counting theorem  with
$$ \left(\alpha_{1}, \alpha_{2}, \alpha_{3}, \alpha_{4}\right)= {\left( -\frac{2(1-H) r}{q}, -\frac{2(1-H) r}{q}, -\frac{2(1-H) (q-r)}{q}, -\frac{2(1-H) (q-r)}{q}\right)}$$
and
$$
(\beta_{1}, \beta_{2}, \beta _{3}, {\beta_{4}})= (-\gamma, -\gamma, -\gamma, -\gamma)$$
with $\gamma>0$ arbitrarly large. We have ($d_{0}$ and $d_{\infty}$ are given by (\ref{d0}) and (\ref{dinf}) respectively)
$$d_{0} (T')= r(T') +\sum _{i=1} ^{4} \alpha _{i} = 3+ 2(2H-2) = 4H-1 >0 \mbox{  for }H>\frac{1}{4}$$
and 
$${d_{\infty} (\emptyset)}= 4-1-4\gamma <0 \mbox{ if }\gamma >\frac{3}{4}.$$
Therefore $T_{1}$ is finite.  Let us regard the last summand, i.e.

\begin{eqnarray*}
T_{2}&:=&\sup_{H\in [\frac{1}{2}, 1]} \int_{\mathbb{R} } dy_{1}....\int_{\mathbb{R} } dy_{4}  {1_{ \vert y_{1} \vert \leqslant \sqrt{2T}} ... 1_{ \vert y_{4} \vert \leqslant\sqrt{2T}}} \\
&&\times \vert y_{1}-y_{2} \vert ^{-\alpha } \vert y_{2}-y_{3} \vert ^{-\alpha } \vert y_{3}-y_{4} \vert ^{-\beta} \vert y_{4}-y_{1} \vert ^{-\beta } <\infty. 
\end{eqnarray*}
This is clearly finite by Lemma 3.3 in \cite{BaiTa} since
$$2\alpha +2\beta +4= 4H-4+4=4H>1$$
when $H>\frac{1}{4}.$

The other summands can be handled by combining the arguments used for the two terms above. For instance, consider

\begin{eqnarray*}
T_{3}&:=&\sup_{H\in [\frac{1}{2}, 1]} \int_{\mathbb{R} } dy_{1}....\int_{\mathbb{R} } dy_{4}  e ^{-\frac{ y_{1}   ^{2}}{2T}}{1_{ \vert y_{1} \vert  \geqslant \sqrt{2T}}  1_{ \vert y_{2} \vert \leqslant \sqrt{2T}}1_{ \vert y_{3} \vert \leqslant \sqrt{2T}} 1_{ \vert y_{4} \vert \leqslant \sqrt{2T}}}\ 
\\
&&\times \vert y_{1}-y_{2} \vert ^{-\alpha } \vert y_{2}-y_{3} \vert ^{-\alpha } \vert y_{3}-y_{4} \vert ^{-\beta} \vert y_{4}-y_{1} \vert ^{-\beta } . 
\end{eqnarray*}
We use the bound (which follows from (\ref{28a-11})

$$y_{1} ^{2}  \geq \frac{1}{4} (\vert y_{1}-y_{2}\vert ^{2}+\vert y_{2}-y_{3}\vert ^{2}+\vert y_{3}-y_{4}\vert ^{2}+\vert y_{4}-y_{1}\vert ^{2}) -(y_{2}^{2}+ y_{3} ^{2}+ y_{4} ^{2})$$
and then

\begin{eqnarray*}
&& e ^{-\frac{y_{1}^{2}}{2T}}\leq e ^{-\frac{1}{8T} (\vert y_{1}-y_{2}\vert ^{2}+\vert y_{2}-y_{3}\vert ^{2}+\vert y_{3}-y_{4}\vert ^{2}+\vert y_{4}-y_{1}\vert ^{2}) } e ^{\frac{y_{2} ^{2}+ y_{3}^{2}+ y_{4} ^{2}}{2T}}\\
&\leq& Ce ^{-\frac{1}{8T} (\vert y_{1}-y_{2}\vert ^{2}+\vert y_{2}-y_{3}\vert ^{2}+\vert y_{3}-y_{4}\vert ^{2}+\vert y_{4}-y_{1}\vert ^{2}) } .
\end{eqnarray*}
The term $T_{3}$  is thus bounded by 

\begin{eqnarray*}
T_{3}&\leq & C\sup_{H\in [\frac{1}{2}, 1]} \int_{\mathbb{R} } dy_{1}....\int_{\mathbb{R} } dy_{4}  e ^{-\frac{1}{8T} (\vert y_{1}-y_{2}\vert ^{2}+\vert y_{2}-y_{3}\vert ^{2}+\vert y_{3}-y_{4}\vert ^{2}+\vert y_{4}-y_{1}\vert ^{2}) } 
\\
&&\times\vert y_{1}-y_{2} \vert ^{-\alpha } \vert y_{2}-y_{3} \vert ^{-\alpha } \vert y_{3}-y_{4} \vert ^{-\beta} \vert y_{4}-y_{1} \vert ^{-\beta }  
\end{eqnarray*}
and we follow the proof for the first term. \qed

\begin{remark}
Notice that the limit  process in Theorem  coincides in distribution with a bifractional Brownian motion with Hurst parameters $H=\frac{1}{2}, K=-\gamma_{0} +1=d-\frac{k}{2}-\sum_{a\in \overline{A}_{k}} H_{a} $ (in the case i. ), $H=\frac{1}{2}, K= d-\frac{k}{2}-\sum_{a\in \overline{A}_{k}} H_{a} + (2H-1)$ (in the case ii.) and $H=K=\frac{1}{2}$ (in the case iii.) We refer to \cite{HaNu}, \cite{T}, \cite{TX} for the definition of the bifractional Brownian motion and for the link between this process and the solution to the heat equation.

\end{remark}

{
\section{Applications to Hermite Ornstein-Uhlenbeck process}
Let $ Z^{q,1}:= Z^{q}$ be a (one-parameter) Hermite process defined by (\ref{hermite-sheet-1}). The Hermite Ornstein Uhlenbeck process has been introduced in \cite{MaTu}.  It is defined as the solution of Langevin equation driven by  Hermite noise. 

\begin{eqnarray}\label{Lang}
X_{t}=\xi -\lambda \int_{0}^{t} X_{s}ds + \sigma Z^{q}_{H}(t) , t \geq  1
\end{eqnarray}
where $\lambda, \sigma > 0 $ and the initial condition $\xi$ is a random variable in $L ^{2}(\Omega)$.
The unique  solution of (\ref{Lang}) is given by

\begin{eqnarray}\label{HOU}
Y^{H}(t)=e^{-\lambda t}\left( \xi + \sigma \int_{0}^{t} e^{ \lambda u} dZ^{q}_{H}(u)\right), \hskip0.5cm t\geq 0
\end{eqnarray}
where the integral $\int_{0}^{t} e^{ \lambda u} dZ^{q}(u)$ exists in the Riemann-Stieljes sense.

In particular, by taking the initial condition $\xi= \sigma \int_{-\infty} ^{0} e ^{\lambda u} d Z ^{H}(u) $ in (\ref{HOU}). The unique solution to (\ref{Lang}), denoted in the sequel by $ (X ^{H}(t))_{t\geq 0}$, can be expressed as
\begin{equation}
\label{HOU2}
X ^{H}(t)= \sigma\int_{-\infty} ^{t}  e ^{-\lambda (t-u) }d Z ^{q}_{H}(u), \hskip0.5cm t\geq 0
\end{equation}
and the stochastic integral in (\ref{HOU2}) can be also understood in the Wiener sense. The process $\left( X ^{H}(t)\right)_{t \geq 0}$ is a stationary process, $H$-self similar process with stationary increments. 

In \cite{SlaTud2} the authors have established the asymptotic behavior with respect to $H$ of the Rosenblatt Ornstein Uhlenbeck process which is the solution of (\ref{Lang}) driven by the Rosenblatt process, i.e.  $q=2$. The proof was based on the analysis of the cumulants, but it is well-known that this method does not work for a Wiener chaos of order $q\geq 3$.  In this section, we will study the behavior as $H\to 1$ and as $H \to \frac{1}{2}$ of the processes $\left( X ^{H}(t)\right)_{ t \in [0,T]}$ and $\left( Y ^{H}(t)\right)_{t \in [0,T]}$ when $q>2$ . The results obtained give a complete picture for the asymptotic behavior of the Hermite Ornstein Uhlenbeck of any order $ q \geq 1$.

\subsection{Asymptotic behavior of the non stationary Hermite Ornstein-Uhlenbeck}
Assume that the initial condition $\xi$ does not depend on $H$.
\begin{prop}\label{pHOU}
\begin{itemize}
\item[1]  Assume $H \to 1$. Then  the process $\left( Y ^{H}(t)\right)_{t \in [0,T]}$ converges weakly, in the space of the continuous functions $C[0,T]$  to the process $\left( Y(t)\right)_{t \in [0,T]}$ given by 
\begin{equation}
Y(t)=e^{-\lambda t}\xi + \sigma\left(1- e ^{-\lambda t }\right)\frac{ H_{q}(Z)}{\sqrt{q!}}
\end{equation}

with $Z \sim \mathcal{N}(0,1)$

\item[2] Assume $H \to \frac{1}{2}$, the process $\left( Y ^{H}(t)\right)_{t \in [0,T]}$ converges weakly, in the space of the continuous functions $C[0,T]$ as $H \to \frac{1}{2}$ to the \textit{standard Ornstein Uhlenbeck process}  $\left( Y_{0}(t)\right)_{t \in [0,T]}$ given by 
\begin{equation}
Y_{0}(t)=e^{-\lambda}\left( \xi+ \sigma \int_{0}^{t} e^{ \lambda u} dW(u)\right)
\end{equation}
that is a  Gaussian process with mean $\mathbf{E} Y_{0}(t)= e ^{-\lambda t} \mathbf{E}\xi$ for any $t\geq 0$ and covariance function
\begin{equation*}
Cov (Y_{0}(t), Y_{0} (s))= \frac{\sigma ^{2}} {2\lambda } \left( e ^{-\lambda \vert t-s\vert}- e ^{-\lambda (t+s)}\right)
\end{equation*}
for every $s,t\geq 0$.

\end{itemize}
\end{prop}
{\bf Proof: } Consider $\alpha_{1}, \ldots , \alpha_{N} \in \mathbb{R}$ and $t_{1}, \ldots , t_{N} \in [0,T].$ We will study the convergence  of the finite dimensional distributions of $Y^{H}$. 

\begin{eqnarray*}
Y_{N}&=&\sum_{i=1}^{N}\alpha_{i}Y^{H}(t_{i})=\sum_{i=1}^{N}e^{-\lambda t_{i}}\xi+ \int_{\mathbb{R}}\sum_{i=1}^{N}\alpha_{i}1_{[0,t_{i}]}(u)e ^{-\lambda (t_{i}-u) }d Z ^{q}_{H}(u)\\
&=& \sum_{i=1}^{N}e^{-\lambda t_{i}}\xi+\int_{\mathbb{R}} f(u)d Z ^{q}_{H}(u)
\end{eqnarray*} 
with $f(u)=\sum_{i=1}^{N}\alpha_{i}1_{[0,t_{i}]}(u)e ^{-\lambda (t_{i}-u) }$.

Notice that in this case the space $\mathcal{H}_{A_{k}}$ given by (\ref{ha}) coincides with $L^{1} (\mathbb{R})$. Since it is clear that $f$  belongs to  
$ \vert \mathcal{H}_{H} \vert \cap L^{1}\left(\mathbb{R}\right)$ (see \cite{SlaTud2}), we get immediatly by Proposition \ref{pp1} the convergence as $H \to 1$ of $\int_{\mathbb{R}} f(u)d Z ^{q}_{H}(u)$ to  $\left(\int_{\mathbb{R}} f(u)du \right) \frac{H_{q}(Z)}{\sqrt{q!}}$.

In order to prove the convergence when $H \to \frac{1}{2}$, we will apply Proposition \ref{pp2}. Using the same arguments as for  the proof of Proposition 5 in \cite{SlaTud2}, we get 
\begin{eqnarray*}
&&\lim _{ H \to \frac{1}{2}  }  H (2 H-1) \int_{\mathbb{R} }\int_{\mathbb{R} }f( u) f( v) \vert { u}-{ v} \vert ^{2{ H}-2}d{ u} d{ v} = \int_{\mathbb{R} } \left(f(u)\right)^{2} du\\
&=& \sum_{i=1}^{N}\sum_{j=1}^{N}\alpha_{i}\alpha_{j}\int_{0 }^{t_{i} \wedge t_{j}}e ^{-\lambda (t_{i}+t_{j}-2u) }du
=\sum_{i=1}^{N}\sum_{j=1}^{N}\alpha_{i}\alpha_{j}\frac{\sigma ^{2}} {2\lambda }\left( e ^{-\lambda \vert t_{i}-t_{j}\vert}- e ^{-\lambda (t_{i}+t_{j})}\right)
\end{eqnarray*} 
which coincides with the variance of $\sum_{j=1} ^{N} \alpha _{j} Y_{0} (t_{j})$. The proof is completed by showing that (\ref{c2-12}) is satisfied. We have 
\begin{eqnarray*}
&&  \int_{\mathbb{R}^{4}} du_{1}...du_{4} f(u_{1})...f(u_{4}) \vert u_{1}- u_{2} \vert ^ {H-1} \vert u_{2}-u_{3} \vert ^ {H-1}\vert u_{3}-u_{4}\vert ^{H-1} \vert u_{4}-u_{1}\vert ^ {H-1}\\
&\leq &\sum_{j_{1},.., j_{4}=1}^ {d} \vert \alpha _{j_{1}}....\alpha_{j_{4}}\vert  \int_{0} ^ {T}...\int_{0} ^ {T} du_{1}..du_{4}  \\
&&\times \vert u_{1}- u_{2} \vert  ^{ \frac{2( H -1)r}{q}} \vert u_{2}-u_{3} \vert  ^{ \frac{2( H -1)r}{q}}\vert u_{3}-u_{4}\vert ^{ \frac{2(H-1)(q-r)}{q}} \vert u_{4}-u_{1}\vert ^{ \frac{2(H-1)(q-r)}{q}}
\end{eqnarray*}
is finite and continuous in $H$ on the set $(\frac{1}{4}, 1]$. This follows from Lemma 3.3 in \cite{BaiTa}  or by applying the power counting theorem with 
$\left(\alpha_{1},\alpha_{2},\alpha_{3}, \alpha_{4}\right)=\left(\frac{2( H -1)r}{q}, \frac{2( H -1)r}{q}, \frac{2(H-1)(q-r)}{q}, \frac{2(H-1)(q-r)}{q}  \right)$.
We recall (see \cite{SlaTud2}) that for $p\geq 1$,
\begin{equation}\label{7d-1}
\mathbf{E}\vert Y^ {H}(t)- Y^ {H}(s) \vert ^ {2p}\leq C_{p} (\mathbf{E} \vert Y^ {H}(t)- Y ^ {H}(s)\vert ^ {2}) ^ {p}\leq c\vert t-s\vert ^ {p}.
\end{equation}
The tighness follows from  (\ref{7d-1}) and Bilingsley criterium (see \cite{B99}).  \qed

\subsection{Asymptotic behavior of the stationary Hermite Ornstein-Uhlenbeck}

Now we will study the asymptotic behavior of (\ref{HOU2}). The diffrence to the non-stationary case is that the function $f$ from the last proof has support of infinite Lebesque measure an we need to use an argument based on the power counting theorem when $H$ tends to one half. The proof of this results is similar in spirit to the proofs  of Proposition 6 and Proposition 7  in \cite{SlaTud2}.
\begin{prop}\label{pHOU2}
\begin{itemize}
\item[1]  Assume $H\to 1$. Then  the process $\left( X ^{H}(t)\right)_{t \in [0,T]}$ converges weakly, in the space of the continuous functions $C[0,T]$  to the process $\left( X(t)\right)_{t \in [0,T]}$ defined by 
\begin{equation}
X(t)=\frac{\sigma}{\lambda}\frac{  H_{q}(Z)}{ \sqrt{q!}}
\end{equation}
with $Z \sim \mathcal{N}(0,1)$
\end{itemize}
\item[2] Assume $H \to \frac{1}{2}$, the process $\left( X ^{H}(t)\right)_{t \in [0,T]}$ converges weakly, in the space of the continuous functions $C[0,T]$ as $H \to \frac{1}{2}$ to the \textit{stationary Ornstein Uhlenbeck process}  $\left( X_{0}(t)\right)_{t \in [0,T]}$ given by 
\begin{equation}
X_{0}(t)= \sigma \int_{-\infty}^{t} e^{- \lambda(t- u)} dW(u)
\end{equation}
which is a  stationary centered Gaussian process with  covariance function
\begin{equation*}
Cov (X_{0}(t), X_{0} (s))= \frac{\sigma ^{2}} {2\lambda }  e ^{-\lambda \vert t-s\vert}
\end{equation*}
for every $s,t\geq 0$.

\end{prop}
{\bf Proof: } Consider $\alpha_{1}, \ldots , \alpha_{N} \in \mathbb{R}$ and $t_{1}, \ldots , t_{N} \in [0,T].$ We will study the convergence  of the finite dimensional distributions of $Y^{H}$. 

\begin{eqnarray*}
\sum_{i=1}^{N}\alpha_{i}X^{H}(t_{i})&=& \int_{\mathbb{R}}\sum_{i=1}^{N}\sigma\alpha_{i}1_{[-\infty,t_{i}]}(u)e ^{-\lambda (t_{i}-u) }d Z ^{q}_{H}(u)\\
&=& \int_{\mathbb{R}} g(u)d Z ^{q}_{H}(u)
\end{eqnarray*} 
with $g(u)=\sum_{i=1}^{N}\alpha_{i}1_{[-\infty,t_{i}]}(u)e ^{-\lambda (t_{i}-u) }$.

The computations in proofs of Proposition 6 and Proposition 7 in \cite{SlaTud2} show that g belongs to  
$ \vert \mathcal{H}_{H} \vert \cap L^{1}\left(\mathbb{R}\right)$, we get immediatly by Proposition \ref{pp1} that the random variable $\sum_{i=1}^{N}\alpha_{i}X^{H}(t_{i}) $  converges to$\sum_{i=1}^{N}\alpha_{i}X(t_{i}) $ as $H \to 1$.

When $H \to \frac{1}{2}$, the proof with slight changes, follows along the same lines as the proof of  Proposition 7 in \cite{SlaTud2}. 
We have 
\begin{equation*}
\mathbf{E} \left( \sum_{j=1} ^ {d}\alpha _{j} X ^ {H}(t_{j})\right) ^{2} \xrightarrow[H \to \frac{1}{2}]{} \mathbf{E} \left( \sum_{j=1} ^ {d} \alpha _{j} X _{0}(t_{j}) \right) ^{2}.
\end{equation*}
It remains to prove that the condition  (\ref{c2-12}) holds true. We have 
\begin{eqnarray*}
&&  \int_{\mathbb{R}^{4}} du_{1}...du_{4} g(u_{1})...g(u_{4}) \vert u_{1}- u_{2} \vert  ^{ \frac{2( H -1)r}{q}} \vert u_{2}-u_{3} \vert  ^{ \frac{2( H -1)r}{q}}\vert u_{3}-u_{4}\vert ^{ \frac{2(H-1)(q-r)}{q}} \vert u_{4}-u_{1}\vert ^{ \frac{2(H-1)(q-r)}{q}}\\
&\leq &\sum_{j_{1}, j_{2},.., j_{4}=1}^ {d} \vert \alpha _{j_{1}}...\alpha _{j_{4}}\vert \int_{-\infty} ^ {t_{j_{1}}} du_{1}....\int_{-\infty} ^ {t_{j_{4}}} du_{m} e ^ {-\lambda (t_{j_{1}}-u_{1})}....e ^ {-\lambda (t_{j_{4}}-u_{4})} \\
&&\vert u_{1}- u_{2} \vert  ^{ \frac{2( H -1)r}{q}} \vert u_{2}-u_{3} \vert  ^{ \frac{2( H -1)r}{q}}\vert u_{3}-u_{4}\vert ^{ \frac{2(H-1)(q-r)}{q}} \vert u_{4}-u_{1}\vert ^{ \frac{2(H-1)(q-r)}{q}}\\
&=&\sum_{j_{1}, j_{2},.., j_{4}=1}^ {d} \vert \alpha _{j_{1}}...\alpha _{j_{4}}\vert \int_{0} ^{\infty} du_{1}....\int_{0}^{\infty} du_{4} e ^ {-\lambda (u_{1}+..+u_{4})}\\
&&\times  \vert u_{1}-u_{2}-(t_{j_{1}}-t_{j_{2}})\vert  ^{ \frac{2( H -1)r}{q}}\vert u_{2}-u_{3}-(t_{j_{1}}-t_{j_{2}})\vert  ^{ \frac{2( H -1)r}{q}} \\
&& \vert u_{3}-u_{4}-(t_{j_{3}}-t_{j_{4}})\vert ^{ \frac{2(H-1)(q-r)}{q}} \vert u_{4}-u_{1}-(t_{j_{4}}-t_{j_{1}})\vert ^{ \frac{2(H-1)(q-r)}{q}}\\
&\leq &e ^{\frac{\lambda }{2}  ( \vert t_{j_{1}}-t_{j_{2}}\vert +...+ \vert t_{j_{4}}-t_{j_{1}}\vert ) }\sum_{j_{1}, j_{2},.., j_{4}=1}^ {d} \vert \alpha _{j_{1}}...\alpha _{j_{4}}\vert \int_{0} ^{\infty} du_{1}...\int_{0} ^{\infty}du_{4} \\
&& e ^{-\frac{\lambda }{2}(\vert u_{1}-u_{2}-(t_{j_{1}}-t_{j_{2}})\vert+...+ \vert u_{4}-u_{1}-(t_{j_{4}}-t_{j_{1}})\vert)}\\
&&\times 
\left( 1\vee   \vert u_{1}-u_{2}-(t_{j_{1}}-t_{j_{2}})\vert  ^{ \frac{2( H -1)r}{q}}\right)\left( 1\vee \vert u_{2}-u_{3}-(t_{j_{1}}-t_{j_{2}})\vert  ^{ \frac{2( H -1)r}{q}} \right)\\
&& \left( 1\vee \vert u_{3}-u_{4}-(t_{j_{3}}-t_{j_{4}})\vert ^{ \frac{2(H-1)(q-r)}{q}} \right) \left( 1\vee \vert u_{4}-u_{1}-(t_{j_{4}}-t_{j_{1}})\vert ^{ \frac{2(H-1)(q-r)}{q}} \right)\\
\end{eqnarray*}

We apply the power counting theorem on the set $T'$ defined by 
$$T' =\{  u_{1}-u_{2}-(t_{j_{1}}-t_{j_{2}}),...,  u_{4}-u_{1}-(t_{j_{4}}-t_{j_{1}}) \}$$ with 
$$(\alpha_{1},..,\alpha _{4})= \left( \frac{2( H -1)r}{q},\frac{2( H -1)r}{q}, \frac{2(H-1)(q-r)}{q}, \frac{2(H-1)(q-r)}{q}\right) \mbox{ and } ( \beta_{1},..,\beta_{4})=(-\gamma, ..., -\gamma)$$
with $\gamma  \in (\frac{3}{4},1]$. Since $T'$ is the only paddet subset of $T'$,  we have
$$d_{0}(T')= 4-1+\frac{4(H-1)(q-r)}{q}+\frac{4(H-1)(q-r)}{q}=4H-1 >0 \mbox{ if } H>\frac{1}{4}$$

\b and 
$$d_{\infty} (\emptyset)=4-1-4\gamma <0 \mbox{ if }\gamma >1-\frac{1}{4}=\frac{3}{4}.$$
Therefore, the function 
$$H\to  \int_{\mathbb{R}}...\int_{\mathbb{R}} du_{1}...du_{4} \vert g (u_{1})...g(u_{m}) \vert \vert u_{1}- u_{2} \vert  ^{ \frac{2( H -1)r}{q}} \vert u_{2}-u_{3} \vert  ^{ \frac{2( H -1)r}{q}}\vert u_{3}-u_{4}\vert ^{ \frac{2(H-1)(q-r)}{q}} \vert u_{4}-u_{1}\vert ^{ \frac{2(H-1)(q-r)}{q}}\\$$
is finite and continuous on the set $D=\{ H \in (0,1],  H>\frac{1}{4}\}$. The conclusion follows from Proposition \ref{pp2}. 

Again the tighness is obtained by (\ref{7d-1}). \qed
}

\section{Appendix}

The basic tools from the analysis on Wiener space and the power counting theorem proven in \cite{TeTa} are presented in this appendix.

\subsection{ Multiple stochastic integrals and the Fourth Moment Theorem}\label{app1}
Here, we shall only recall some elementary
facts; our main reference is  \cite{N}. Consider
${\mathcal{H}}$ a real separable infinite-dimensional Hilbert space
with its associated inner product ${\langle
.,.\rangle}_{\mathcal{H}}$, and $(B (\varphi),
\varphi\in{\mathcal{H}})$ an isonormal Gaussian process on a
probability space $(\Omega, {\mathfrak{F}}, \mathbb{P})$, which is a
centered Gaussian family of random variables such that
$\mathbf{E}\left( B(\varphi) B(\psi) \right) = {\langle\varphi,
\psi\rangle}_{{\mathcal{H}}}$, for every
$\varphi,\psi\in{\mathcal{H}}$. Denote by $I_{q}$ the $q$th multiple
stochastic integral with respect to $B$. This $I_{q}$ is actually an
isometry between the Hilbert space ${\mathcal{H}}^{\odot q}$
(symmetric tensor product) equipped with the scaled norm
$\frac{1}{\sqrt{q!}}\Vert\cdot\Vert_{{\mathcal{H}}^{\otimes q}}$ and
the Wiener chaos of order $q$, which is defined as the closed linear
span of the random variables $H_{q}(B(\varphi))$ where
$\varphi\in{\mathcal{H}},\;\Vert\varphi\Vert_{{\mathcal{H}}}=1$ and
$H_{q}$ is the Hermite polynomial of degree $q\geq 1$ defined
by:\begin{equation}\label{Hermite-poly}
H_{q}(x)=(-1)^{q} \exp \left( \frac{x^{2}}{2} \right) \frac{{\mathrm{d}}^{q}%
}{{\mathrm{d}x}^{q}}\left( \exp \left(
-\frac{x^{2}}{2}\right)\right),\;x\in \mathbb{R}.
\end{equation}The isometry of multiple integrals can be written as: for $p,\;q\geq
1$,\;$f\in{{\mathcal{H}}^{\otimes p}}$ and
$g\in{{\mathcal{H}}^{\otimes q}}$,
\begin{equation} \mathbf{E}\Big(I_{p}(f) I_{q}(g) \Big)= \left\{
\begin{array}{rcl}\label{iso}
q! \langle \tilde{f},\tilde{g}
\rangle _{{\mathcal{H}}^{\otimes q}}&&\mbox{if}\;p=q\\
\noalign{\vskip 2mm} 0 \quad\quad&&\mbox{otherwise}.
\end{array}\right.
\end{equation}It also holds that:
\begin{equation*}
I_{q}(f) = I_{q}\big( \tilde{f}\big),
\end{equation*}
where $\tilde{f} $ denotes the canonical symmetrization of $f$ and it is defined by: $$\tilde{f}%
(x_{1}, \ldots , x_{q}) =\frac{1}{q!}\sum_{\sigma\in\mathcal{S}_q}
f(x_{\sigma (1) },\ldots, x_{\sigma (q)}),$$in which the sum runs
over all permutations $\sigma$ of $\{1,\ldots,q\}$.

In the particular case when $\mathcal{H}=L^2(T, \mathcal{B}(T), \mu)$ , the $r$th
contraction $f\otimes_{r}g$ is the element of
${\mathcal{H}}^{\otimes(p+q-2r)}$, which is defined by: {
\begin{eqnarray}\label{contra}
& (f\otimes_{r} g) ( s_{1}, \ldots, s_{p-r}, t_{1}, \ldots, t_{q-r})  \notag \\
& =\int_{T ^{r} }\mathrm{d}u_{1}\ldots \mathrm{d}u_{r} f( s_{1},
\ldots, s_{p-r}, u_{1}, \ldots,u_{r})g(t_{1}, \ldots, t_{q-r},u_{1},
\ldots,u_{r}),
\end{eqnarray}for every $f\in L^2({[0,T]}^p)$, $g\in L^2({[0,T]}^q)$
and $r=1,\ldots,p\wedge q$.

An important  property of  finite sums of multiple integrals is the hypercontractivity. Namely, if $F= \sum_{k=0} ^{n} I_{k}(f_{k}) $ with $f_{k}\in \mathcal{H} ^{\otimes k}$ then
\begin{equation}
\label{hyper}
\mathbf{E}\vert F \vert ^{p} \leq C_{p} \left( \mathbf{E}F ^{2} \right) ^{\frac{p}{2}}.
\end{equation}
for every $p\geq 2$.

We will use the following  famous result initially proven in \cite{NuPe} that characterizes the convergence in distribution of a sequence of multiple integrals torward the Gaussian law.

\begin{theorem}\label{fmt}
Fix $n \geq 2$ and let $\left( F_{k} , k \geq 1 \right)$ , $F_{k} = I_{n}\left(f_{k} \right) $ ( with $f_{k} \in {\mathcal{H}}^{\odot n}$ for every $k \geq 1$ ), be a sequence of square-integrable random variables in the nth Wiener chaos such that $ \mathbf{E} \left[ F_{k}^{2} \right] \rightarrow 1$ as $k \rightarrow \infty$. The following are equivalent:
\begin{enumerate}
\item the sequence $\left( F_{k} \right)_{k \geq 0}$ converges in distribution to the normal law $\mathcal{N} (0,1)$; 
\item $ \mathbf{E} \left[ F_{k}^{4} \right] = 3$ as $k \rightarrow \infty$;
\item for all  $1 \leq  l \leq n-1$, it holds that $\lim\limits_{k \rightarrow \infty} \Vert  f_{k}  \otimes_{l}  f_{k} \Vert_{{\mathcal{H}}^{\otimes 2 (n-l)}} = 0 $;

\end{enumerate}
\end{theorem}
Another equivalent condition can be stated in term of the Malliavin derivatives of $F_{k}$, see \cite{NPbook}.

\subsection{Power counting theorem}

We need to recall some notation and results from \cite{TeTa} which are needed in order to check the integrability assumption from Proposition  \ref{pp2}.

Consider a set $T=\{M_{1},.., M_{m}\}$ of linear functions on $\mathbb{R } ^ {m}$.  The power counting theorem (see Theorem 1.1 and Corollary 1.1 in \cite{TeTa}) gives sufficient conditions for the integral 
\begin{equation}
\label{i}
I= \int_{\mathbb{R}}...\int_{\mathbb{R}} du_{1}...du_{m} f_{1} (M_{1} (u_{1},.., u_{m}))....f_{m} (M_{m} (u_{1},.., u_{m}))
\end{equation}
to be finite, where $f_{i}:\mathbb{R} \to \mathbb{R}$, $i=1,.., m$ are such that $\vert f_{i} \vert $ is bounded above on $(a_{i}, b_{i})$ ($0<a_{i}<b_{i} <\infty$) and
$$ \vert f_{i}(y) \vert \leq c_{i} \vert y\vert ^ {\alpha _{i} } \mbox{ if } \vert y_{i}\vert <a_{i} \mbox{ and } \vert f_{i}(y)\vert \leq c_{i}\vert y\vert ^ {\beta _{i} } \mbox{ if } \vert y\vert >b_{i}.$$

For a subset $W\subset T$ we denote  by $s_{T}(W) =span (W)\cap T$. A subset $W$ of $T$ is said to be {\em padded } if $s_{T}(W)=W$ and any functional $M\in W$ also belongs to $s_{T}(W\setminus \{M\}).$ Denote by $\span (W)$ the linear span generated by $W$ and by $r(W)$ the number of linearly independent elements of $W$. 

Then Theorem 1.1  in \cite{TeTa} says that the integral $I$ (\ref{i}) is finite if
\begin{equation}
\label{d0}
 d_{0} (W)= r(W) + \sum _{s_{T}(W)} \alpha _{i} >0 
\end{equation}
for any  subset $W$ of $T$ with   $s_{T}(W)=W$ and 
\begin{equation}
\label{dinf}
 d_{\infty } (W)= r(T) -r(W) + \sum _{T\setminus s_{T}(W)} \beta _{i} <0
\end{equation}
for any proper subset $W$ of $T$ with $s_{T}(W)= W$, including the empty set. If $\alpha _{i}>-1$ then  it suffices to check (\ref{d0}) for any padded subset $W\subset T$. Also, it suffices to verify (\ref{dinf}) only for padded subsets of $T$ if $ \beta _{i} \geq -1.$

The condition (\ref{d0}) implies the integrability at the origin while (\ref{dinf}) gives the integrability of $I$ at infinity.

There is a similar result if one starts with a set $T$ of affine functionals instead of linear functionals.

\end{document}